\documentclass[reqno]{amsart}
\usepackage{amsmath}%
\usepackage{amsfonts}%
\usepackage{amssymb}%
\usepackage{graphicx}
\usepackage{setspace}
\usepackage{microtype}

\theoremstyle{plain}

\numberwithin{equation}{section}
\begin{document}
\onehalfspacing
\raggedbottom
\title[Fact Sheet, Part I]{Fact Sheet Research on Bayesian Decision Theory}
\author{H.R.N.~van~Erp}
\author{R.O.~Linger}


\author{P.H.A.J.M.~van~Gelder}


\begin{abstract}
In this fact sheet we give some preliminary research results on the Bayesian Decision Theory. This theory has been under construction for the past two years. But what started as an intuitive enough idea, now seems to have the makings of something more fundamental.  
\end{abstract}

\maketitle

\section{Introduction}
It has been shown that the product and sum rules of both the Bayesian probability and the Bayesian information theories are derivable by way of consistency, \cite{Cox46, Jaynes03, vanErp13, Knuth10, Skilling08}; the implication being that in our plausibility and relevancy judgments we humans have a preference for consistency, or, equivalently, rationality. Moreover, Knuth, a MaxEnt-Bayesian\footnote{MaxEnt-Bayesians are, as a rule, physicists that trace their statistical lineage from Jaynes, back to Jeffreys, back to Laplace.}, is now researching if the very laws of physics may be derived by way of consistency constraints on lattices of events, \cite{Knuth14}. 

In light of these exciting new developments, in both the fields of inference and physics, where consistency arguments are taking the stage by storm, and the fact that these authors, after two years of continuous research, have reached the point that they have come to trust their Bayesian decision theory to almost the same extent as they have grown to trust the Bayesian probability and information theories\footnote{The former being their field of expertise, and the latter being the subject matter of the first author's current thesis work.}, these authors have come to entertain the notion that maybe their Bayesian decision theory, which initially started as an intuitive enough Bayesian alternative for the paradigm of behavioral economics, might actually be Bayesian in the strictest sense of the word; that is, an inescapable consequence of the desideratum of consistency.

In this fact sheet we will present the case, as it currently stands, for the Bayesian decision theory, and which led us to this conjecture, together with the consistency proof of the Bernoulli utility function, also known as the Weber-Fechner law of sense perception.  


\section{The Bayesian Decision Theory}
The Bayesian decision theory is very simple in structure. Its algorithmic steps are the following:
\begin{enumerate}
	\item Use the product and sum rules of Bayesian probability theory to construct outcome probability distributions. \\
	\item If our outcomes are monetary in nature, then by way of the Bernoulli utility function we may map utilities\footnote{Or, as Bernoulli called them, moral values, \cite{Bernoulli38}.} to the monetary outcomes of our outcome probability distributions. \\
	\item Maximize the sum of the lower and upper bounds of the resulting utility probability distributions. 
\end{enumerate}
This, then, is the whole of the Bayesian decision theory.

\section{Constructing Outcome Probability Distributions}
\label{chapter2.1}
In the Bayesian decision theory each problem of choice is understood to consist of a set of decisions from which we must choose. Each possible decision, when taken, has its own set of possible outcomes, and each outcome, for a given decision, has its own plausibility of occurring relative to the other outcomes under that same decision. So, each decision in our problem of choice admits its own outcome probability distribution.

We will demonstrate in this section how to construct outcome probability distributions using the rules of Bayesian probability theory\footnote{See also Appendix~\ref{Bayes}.}. In our hypothetical problem of choice, the possible decisions $D_{i}$, for $i = 1, 2$, under consideration are whether or not to wear seat belts:
\[
	D_{1} = \text{Seat belts}, \qquad		D_{2} = \text{No Seat belts}.
\]
The relevant events $E_{j}$, for $j = 1, 2, 3$, when driving a car, as perceived by the decision maker, are
\[
	E_{1} = \text{No Accident}, \quad		E_{2} = \text{Small Accident}, \quad		E_{3} = \text{Severe Accident}.
\]
The perceived outcomes $O_{k}$, for $k = 1, 2, 3$, are
\[
	O_{1} = \text{No Harm}, \quad		O_{2} = \text{Some Bruises}, \quad		O_{3} = \text{Broken Bones}.
\]
For this particular case, the decisions taken do not modulate the probabilities of an event. So, we have that the probability for an event conditional on the decision taken is the same for both decisions, say:
\begin{equation}
\label{eq.P1.2.1}
	P\!\left(\left.E_{1}\right|D_{i}\right) = 0.950, \quad		P\!\left(\left.E_{2}\right|D_{i}\right) = 0.049, \quad		P\!\left(\left.E_{3}\right|D_{i}\right) = 0.001,
\end{equation}
for $i = 1, 2$. However, the conditional probability distributions of the outcomes given an event are modulated by the decision taken. 

We first consider the case where the decision maker is considering to wear seat belts, that is, $D_{1}$. Say, we have the following conditional probabilities: 
\begin{align}
	\label{eq.P1.2.2}
&P\!\left(\left.O_{1}\right|E_{1}, D_{1}\right) = 1.00,  	& &P\!\left(\left.O_{2}\right|E_{1}, D_{1}\right) = 0.00, 	& 	&P\!\left(\left.O_{3}\right|E_{1}, D_{1}\right) = 0.00, \nonumber \\
&P\!\left(\left.O_{1}\right|E_{2}, D_{1}\right) = 0.75, 	& &P\!\left(\left.O_{2}\right|E_{2}, D_{1}\right) = 0.25, 	&	  &P\!\left(\left.O_{3}\right|E_{2}, D_{1}\right) = 0.00, \nonumber \\
&P\!\left(\left.O_{1}\right|E_{3}, D_{1}\right) = 0.20, 	& &P\!\left(\left.O_{2}\right|E_{3}, D_{1}\right) = 0.70, 	& 	&P\!\left(\left.O_{3}\right|E_{3}, D_{1}\right) = 0.10. \nonumber \\
\end{align}
Then by way of the product rule, \cite{Jaynes03},
\begin{equation}
	\label{eq.P1.2.3}
	 P\!\left(A\right) P\!\left(\left.B\right|A\right) = P\!\left(A B\right) = P\!\left(B\right) P\!\left(\left.A\right|B\right),
\end{equation}
we may combine the probability of an event, \eqref{eq.P1.2.1}, with the corresponding conditional probability distributions of some outcome given that event, \eqref{eq.P1.2.2}, and obtain the probabilities of an event $E_{j}$ and an outcome $O_{k}$ given decision $D_{i}$: 
\begin{equation}
	\label{eq.P1.2.3b}
	P\!\left(\left.E_{j}, O_{k}\right|D_{i}\right) = P\!\left(\left.E_{j}\right|D_{i}\right) P\!\left(\left.O_{k}\right|E_{j}, D_{i}\right).
\end{equation}
We may present all these probabilities \eqref{eq.P1.2.3b} in a table and so get the corresponding bivariate probability distribution, Table~\ref{tab:biDist1}.
\begin{table}[h]
	\centering
		\begin{tabular}{c|c|c|c}
		
						&															&													& 								\\
$P\!\left(\left.E_{j}, O_{k}\right| D_{1}\right)$			&	$O_{1} = \text{No Harm}$						&	$O_{2} = \text{Some Bruises}$								& $O_{3} = \text{Broken Bones}$				\\
			 			& 														& 							  				& 								\\ \hline
			 			& 														& 							  				& 								\\
$E_{1} = \text{No Accident}$:		& 0.9500 											& 0.0											& 0.0							\\
 						& 														& 							  				& 								\\ \hline
 						& 														&  												&									\\ 
$E_{2} = \text{Small Accident}$:		& 0.0370  										& 0.0120									& 0.0		 					\\ 
						& 	 													&   											&									\\ \hline
						& 														& 							  				& 								\\
$E_{3}  = \text{Severe Accident}$:		& 0.0002 											& 0.0007									& 0.0001					\\ 
  					&  														&  				   							&									\\ \hline \\
  					
		\end{tabular}	
	\caption{Bivariate event-outcome probability distribution for $D_{1}$ }
	\label{tab:biDist1}
\end{table}

Let $A_{i} = \left\{A_{1},\ldots,A_{n}\right\}$  be a set of $n$ mutually exclusive and exhaustive propositions, that is, one and only one of the $A_{i}$ is necessarily true. Let $B_{j} = \left\{B_{1},\ldots,B_{m}\right\}$  be another set of $m$ mutually exclusive and exhaustive propositions. Then, by way of the generalized sum rule, \cite{Jaynes03}, we have
\begin{equation}
	\label{eq.P1.2.4}
	\sum_{i=1}^{n} P\!\left(A_{i},B_{j}\right) = P\!\left(A_{1},B_{j}\right) + \cdots + P\!\left(A_{n},B_{j}\right) = P\!\left(B_{j}\right),
\end{equation}
where $\sum_{j} P\!\left(B_{j}\right) = 1$. Using this generalized sum rule, we may `marginalize' the event-outcome probabilities $P\!\left(\left.E_{j}, O_{k}\right| D_{1}\right)$ over the events $E_{j}$, that is,
\begin{equation}
	\label{eq.P1.2.4b}
	P\!\left(\left.O_{k}\right| D_{i}\right) = \sum_{j=1}^{m} P\!\left(\left.E_{j}, O_{k}\right| D_{i}\right) ,
\end{equation}
and so get the marginalized outcome probability distribution, Table~\ref{tab:uniDist1}. 
\begin{table}[h]
	\centering
		\begin{tabular}{c|c|c|c}
						
						&															&													& 								\\
						&	$O_{1} = \text{No Harm}$						&	$O_{2} = \text{Some Bruises}$								& $O_{3} = \text{Broken Bones}$				\\
			 			& 														& 							  				& 								\\ \hline
			 			& 														& 							  				& 								\\
$P\!\left(\left. O_{k}\right| D_{1}\right)$		& 0.9872										& 0.0127								& 0.0001				\\
 						& 														& 							  				& 								\\ \hline \\

		\end{tabular}
	\caption{Marginalized outcome probability distribution for $D_{1}$ }
	\label{tab:uniDist1}
\end{table}

We now consider the case we the decision maker is considering not to wear seat belts, that is, $D_{2}$. Say we have the following conditional probabilities: 
\begin{align}
	\label{eq.P1.2.5}
&P\!\left(\left.O_{1}\right|E_{1}, D_{2}\right) = 1.00,  	&  &P\!\left(\left.O_{2}\right|E_{1}, D_{2}\right) = 0.00, 	& 	&P\!\left(\left.O_{3}\right|E_{1}, D_{2}\right) = 0.00, \nonumber \\
&P\!\left(\left.O_{1}\right|E_{2}, D_{2}\right) = 0.25, 	&  &P\!\left(\left.O_{2}\right|E_{2}, D_{2}\right) = 0.75, 	& 	  &P\!\left(\left.O_{3}\right|E_{2}, D_{2}\right) = 0.00, \nonumber \\
&P\!\left(\left.O_{1}\right|E_{3}, D_{2}\right) = 0.10, 	&  &P\!\left(\left.O_{2}\right|E_{3}, D_{2}\right) = 0.30, 	&  	&P\!\left(\left.O_{3}\right|E_{3}, D_{2}\right) = 0.60. \nonumber \\
\end{align}
Then, using \eqref{eq.P1.2.3}, we may combine the probability of an event, \eqref{eq.P1.2.1}, with the corresponding conditional probability distributions of some outcome given that event, \eqref{eq.P1.2.5}, and so get the corresponding bivariate probability distribution, Table~\ref{tab:biDist2}.
\begin{table}[h]
	\centering
		\begin{tabular}{c|c|c|c}
		
						&															&													& 								\\
$P\!\left(\left.E_{j}, O_{k}\right| D_{2}\right)$			&	$O_{1} = \text{No Harm}$						&	$O_{2} = \text{Some Bruises}$								& $O_{3} = \text{Broken Bones}$				\\
			 			& 														& 							  				& 								\\ \hline
			 			& 														& 							  				& 								\\
$E_{1} = \text{No Accident}$:		& 0.9500 											& 0.0											& 0.0							\\
 						& 														& 							  				& 								\\ \hline
 						& 														&  												&									\\ 
$E_{2} = \text{Small Accident}$:		& 0.0120 										& 0.0370									& 0.0		 					\\ 
						& 	 													&   											&									\\ \hline
						& 														& 							  				& 								\\
$E_{3}  = \text{Severe Accident}$:		& 0.0001 											& 0.0003									& 0.0006					\\ 
  					&  														&  				   							&									\\ \hline \\
  				
		\end{tabular}
	\caption{Bivariate event-outcome probability distribution for $D_{2}$ }
	\label{tab:biDist2}
\end{table}

\noindent Marginalizing the event-outcome probabilities $P\!\left(\left.E_{j}, O_{k}\right| D_{2}\right)$ over the events $E_{j}$, \eqref{eq.P1.2.4}, we get the marginalized outcome probability distribution, Table~\ref{tab:uniDist2}.
\begin{table}[h]
	\centering
		\begin{tabular}{c|c|c|c}
						
						&															&													& 								\\
						&	$O_{1} = \text{No Harm}$						&	$O_{2} = \text{Some Bruises}$								& $O_{3} = \text{Broken Bones}$				\\
			 			& 														& 							  				& 								\\ \hline
			 			& 														& 							  				& 								\\
$P\!\left(\left. O_{k}\right| D_{1}\right)$		& 0.9621										& 0.0373								& 0.0006				\\
 						& 														& 							  				& 								\\ \hline \\
 					
		\end{tabular}
	\caption{Marginalized outcome probability distribution for $D_{2}$ }
	\label{tab:uniDist2}
\end{table}

In its most abstract form, we have that each problem of choice consists of a set of potential decisions
\[
	D_{i} = \left\{D_{1},\ldots,D_{n} \right\}.
\] 
Each decision $D_{i}$ we make may give rise to a set of possible events 
\[
	E_{j_{i}} = \left\{E_{1_{i}},\ldots, E_{m_{i}} \right\}.
\]
These events $E_{j_{i}}$ are associated with the decisions $D_{i}$ by way of the conditional probabilities $P\!\left(\left.E_{j_{i}}\right|D_{i}\right)$. Furthermore, each event $E_{j_{i}}$ allows for a set of potential outcomes
\[
	O_{k_{j_{i}}} = \left\{O_{1_{j_{i}}},\ldots, O_{l_{j_{i}}}\right\}. 
\]
These outcomes $O_{k_{j_{i}}}$ are associated with the events $E_{j_{i}}$ by way of the conditional probabilities $P\!\left(\left.O_{k_{j_{i}}}\right|E_{j_{i}}\right)$. 

By way of the product rule, \eqref{eq.P1.2.3}, we compute the bivariate probability distribution of an event and an outcome conditional on the decision taken:
\begin{equation}
	\label{eq.P1.2.6}
	P\!\left(\left.E_{j_{i}}, O_{k_{j_{i}}}\right| D_{i}\right)= P\!\left(\left.E_{j_{i}}\right|D_{i}\right) P\!\left(\left.O_{k_{j_{i}}}\right|E_{j_{i}}\right).
\end{equation}
The outcome probability distribution is then obtained by marginalizing, \eqref{eq.P1.2.4}, over all the possible events
\begin{equation}
	\label{eq.P1.2.7}
	P\!\left(\left. O_{k_{j_{i}}}\right| D_{i}\right) = \sum_{j_{i} = 1}^{m_{i}} P\!\left(\left.E_{j_{i}}, O_{k_{j_{i}}}\right| D_{i}\right).
\end{equation}
The outcome probability distributions \eqref{eq.P1.2.7}, for $i = 1,\ldots, n$, are the information carriers which represent our state of knowledge in regards to the consequences of our decisions.

At first sight the added event space in the abstract form given here may seem somewhat superfluous\footnote{As this event space hardly is used in this fact sheet.}. But it was felt that further down the line, in decision theoretical problems more complex than the ones given in this fact sheet, this added space may help one in the construction of outcome probability distributions\footnote{If we want to collect all the possible outcomes in one probability distribution, under a given decision, then we may first gather the different events, which precede these outcomes, in a seperate probability distribution. We foresee that this factoring may greatly help in the construction of complex outcome probability distributions; cascading events, event trees, etc...}.

\section{Translating Monetary Outcomes to Utilities}
The Bernoulli utility function in the field of psycho-physics\footnote{Psycho-physics is the experimental field of psychology that studies sense perception.} is called Weber-Fechner law. Seeing that we will use in this section the psycho-physical point of view of money increments as a stimulus, we will refer in what follows to the Bernoulli utility function as the Weber-Fechner law. But both names point to the same function\footnote{In the next section we will formally derive the Bernoulli utility function by way of a novel consistency argument. In Appendices~\ref{BernoulliA} and~\ref{BernoulliB} we, respectively, discuss the ubiquitousness of the Bernoulli utility function, which has been derived several times over by different arguments, and proceed to give its corollary, the negative Bernoulli utility function, which may be used to model the utility of debt.}.

The translation of monetary stimuli to utilities is analogous to the case where we are asked to translate loudness to a numerical value. According to Weber-Fechner law, postulated in the 19th century\footnote{Which is just Bernoulli's utility function, postulated in the 18th century; see Appendix~\ref{BernoulliA}.} by the experimental psychologist Fechner, intuitive human sensations tend to be logarithmic functions of the difference in stimulus, \cite{Fechner60}. So, we do not perceive stimuli in isolation, rather we perceive the relative change in stimuli, case in point being the decibel scale of sound. 

Let $S_{1}$ and $S_{2}$ be two stimuli which are to be compared. Then the Weber-Fechner law tells us that the Relative Change (RC) is the difference of the logarithms of the stimuli:
\begin{equation}
	\label{eq.P1.5.1}
	\text{RC} = c \log_{d} S_{2} - c \log_{d} S_{1} = c \log_{d} \frac{S_{2}}{S_{1}},
\end{equation}
where  $c$ is some scaling factor and $d$ some base of the logarithm. From \eqref{eq.P1.5.1}, we have that if stimuli $S_{1}$ and $S_{2}$ are indistinguishable, that is, of the same strength, then their RC is 0. If $S_{2}$ increases relative to $S_{1}$, then $\text{RC} > 0$. If $S_{2}$ decreases relative to $S_{1}$, then $\text{RC} < 0$.

The Weber-Fechner law allows for one degree of freedom. This can be seen as follows. Since 
\[
	\log_{d} x =  \frac{\log x}{\log d},
\]
we can rewrite \eqref{eq.P1.5.1} as 
\begin{equation}
	\label{eq.P1.5.2}
	\text{RC} = q \log \frac{S_{2}}{S_{1}},
\end{equation}
where
\begin{equation}
	\label{eq.P1.5.3}
	q = \frac{c}{\log d}.
\end{equation}

Let $\Delta S$ be an increment, either positive or negative, in a monetary stimulus $S$. Then we may define the utility of a monetary increment $\Delta S$ to be the perceived relative change in the initial wealth $S$ due to that increment $\Delta S$, \eqref{eq.P1.5.2}:
\begin{equation}
	\label{eq.P1.5.3b}
	u\!\left(\left.\Delta S\right|S\right) = q \log \frac{S + \Delta S}{S} , \qquad			 -S < \Delta S < \infty.
\end{equation}

If $\Delta S = -S$, then \eqref{eq.P1.5.3b} tells us that a loss of all one's initial wealth $S$ would have a utility of minus infinity. This is clearly not realistic. So, in order to model such a loss, we must introduce the threshold of income which is still significant $\gamma$, \cite{Jaynes03}, where $\gamma > 0$. The threshold of income has the following interpretation. 

Even for the homeless person there is some minimum amount of money that is still significant. This may be one dollar for a bag of potato chips, or three dollars for a packet of cigarettes. If the loss of money breaks through the limit of the minimum significant amount $\gamma$, the homeless person is left with an amount of money which, for all intents and purposes, is worthless. Using the concept of the threshold of income, we may modify \eqref{eq.P1.5.3b} as 
\begin{equation}
	\label{eq.P1.5.4b}
	u\!\left(\left.\Delta S\right|S\right) = q \log \frac{S + \Delta S}{S} , \qquad			 -S + \gamma < \Delta S < \infty.
\end{equation} 

If we want to give a graphical representation of \eqref{eq.P1.5.4b}, then the scaling constant $q$, also known as the Weber constant, must be set to some numerical value. 

Say, we have a monthly expendable income of a thousand dollars, for groceries and the like, then introspection\footnote{Introspection being the starting point of all psychological experimentation.} would suggest that a loss or gain of an amount less than ten dollars would not move us that much. 

So, $\Delta S = 10$ constitutes a just noticeable difference, or, equivalently, 1 utile, for an initial wealth of $S =1000$, \eqref{eq.P1.5.4b}:
\begin{equation}
	\label{eq.P1.10.2a}
	1 \:\text{utile} = q \log \frac{1000 + 10}{1000}.
\end{equation} 
If we then solve for the unknown Weber constant $q$, we find 
\begin{equation}
	\label{eq.P1.10.2b}
		q =\frac{1}{\log 1010 - \log 1000}\approx 100.
\end{equation} 
Note that utiles represent the utility of the monetary outcomes, much like decibels represent the perceived intensity of sound\footnote{Note that for the decibel scale the Weber constant has been determined to be $q = 10/\log 10 = 4.34$.}.

\begin{figure}[h]
		\centering
			\includegraphics[width=0.40\textwidth]{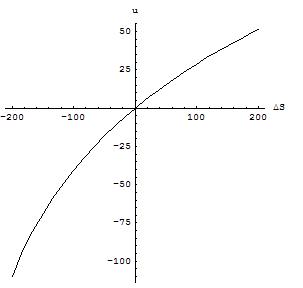}
		\caption{Utility plot for initial wealth 300 dollars}
		\label{fig:P1.5.1}
	\end{figure}

Suppose we have a student who has three hundred dollars per month to spend on groceries and the like and who stands to lose or to gain up to two hundred dollars. Then, by way of \eqref{eq.P1.5.4b} and \eqref{eq.P1.10.2b}, we obtain the following mapping of monetary outcomes to utilities, Figure~\ref{fig:P1.5.1}.

\begin{figure}[h]
		\centering
			\includegraphics[width=0.40\textwidth]{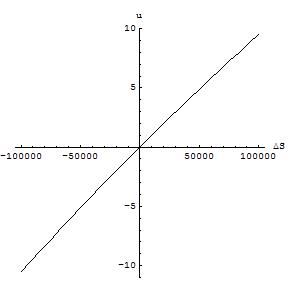}
		\caption{Utility plot for initial wealth 1.000.000 dollars}
		\label{fig:P1.5.2}
	\end{figure}

\noindent For the case of the rich man who has one million dollars to spend on groceries and the like and who stands stands to lose or to gain up to a hundred thousand dollars, we obtain the alternative mapping, Figure~\ref{fig:P1.5.2}.

Loss aversion is the phenomenon that losses may loom larger than gains, \cite{Tversky92}. Comparing Figures~\ref{fig:P1.5.1} and~\ref{fig:P1.5.2}, we see that the Weber-Fechner law of experimental psychology, which is just Bernoulli's utility function, captures both the loss aversion of the poor student, that is, asymmetry in gains and losses, as well as the linearity of the utility of relatively small gains and losses for the rich man. 

\section{A Consistency Proof of the Bernoulli Utility Function}
\label{BernoulliC}
We will derive the Bernoulli utility function, or, equivalently, the Weber-Fechner law, using the desiderata of invariance and consistency. In this we follow a venerable Bayesian tradition, \cite{Cox46, Jaynes03, Knuth10}.

\subsection{Setting up the functional equations}
Say, we have the positive quantities $x$, $y$, $z$, $t$ of some stimulus or commodity of interest. These quantities, being numbers on the positive real and, consequently, transitive, admit an ordering. The quantities are assumed to be ordered as $x \leq y \leq z \leq t$. We now want to find the function $f$ that quantifies the perceived decrease associated with going from, say, the quantity $y$ to the quantity $x$.

The first functional equation is based on the observation that the unknown function $f$ should be invariant for a change in the unit of the quantity, that is,
\begin{equation}
	\label{BerC.1}
	f\!\left(c x, c y\right) = f\!\left(x, y\right),
\end{equation}
where $c$ is some scaling factor. 

For, example, if our quantities concern sums of money, then a perceived loss should be invariant for a change in the monetary unit. The perceived loss of going from ten dollars to one dollar should be the same perceived loss if we reformulate this scenario in dollar cents, that is, if we go from a thousand dollar cents to a hundred dollar cents.

Note that we also may interpret \eqref{BerC.1} in the following way. If person $A$ has an initial wealth of $y$ and person $B$ has a wealth of $c y$, where $c$ is some positive constant, then the perceived loss of person $A$ when going from $y$ to $x$ should be the same as the perceived loss of person $B$ when going from $c y$ to $c x$. This alternative interpretation of \eqref{BerC.1} connects with Bernoulli's original derivation of the utility function, in that it is a corollary of the reasoning that led him to his third and final consideration\footnote{See Appendix~\ref{BernoulliA}.}.

In order to give the second functional equation, we must introduce the function $g$ that enforces a change of context, \cite{Knuth10}:
\begin{equation}
	\label{BerC.2}
	f\!\left(x, z\right) = g\!\left[f\!\left(x, y\right), f\!\left(y, z\right)\right].
\end{equation}
The functional equation \eqref{BerC.2} is a consistency equation that states that the perceived decrease in going from $z$ to $x$, should be the same perceived decrease if we first go from $z$ to $y$, and then from $y$ to $x$.

If we assume differentiability, then the invariance equations \eqref{BerC.1} and \eqref{BerC.2}, together with the two boundary conditions:
\begin{equation}
	\label{BerC.3}
	f\!\left(x, x\right) = 0,
\end{equation}
and
\begin{equation}
	\label{BerC.4}
	f\!\left(x, y\right) < 0,  \quad \text{for } x < y,
\end{equation}
are sufficient to find the function $f$ that quantifies the perceived decrease associated with going from the quantity $y$ to the quantity $x$.

\subsection{Solving the functional equations}
We first will look at the structure that the invariance equation \eqref{BerC.1} provides us in determining the form of $f$. To do so, we take for \eqref{BerC.1} the derivative over $c$:
\begin{equation}
	\label{BerC.5}
	\frac{d}{dc} f\!\left(c x, c y\right) = \frac{d}{dc}  f\!\left(x, y\right)
\end{equation}
which leaves with the partial differential equation
\begin{equation}
	\label{BerC.6}
	x f^{\left(1, 0\right)}\!\left(c x, c y\right) + y f^{\left(0, 1\right)}\!\left(c x, c y\right) = 0
\end{equation}
or, equivalently, if we let $u = c x$ and $v = c y$,
\begin{equation}
	\label{BerC.7}
	\frac{u}{c} f^{\left(1, 0\right)}\!\left(u, v\right) + \frac{v}{c} f^{\left(0, 1\right)}\!\left(u, v\right) = 0.
\end{equation}
This partial differential equation has as its general solution:
\begin{equation}
	\label{BerC.8}
		f\!\left(u, v\right) = h\!\left(\frac{u}{v}\right),
\end{equation}
were $h$ is some arbitrary function.

The consistency equation \eqref{BerC.2} may be solved as follows. First, by way of \eqref{BerC.2}, we construct the associativity equation, \cite{Aczel66, Knuth10}:
\begin{equation}
	\label{BerC.9}
		g\!\left\{g\!\left[f\!\left(x, y\right), f\!\left(y, z\right)\right], f\!\left(z, t\right)\right\} = f\!\left(x, t\right) = g\!\left\{f\!\left(x, y\right), g\!\left[f\!\left(y, z\right), f\!\left(z, t\right)\right]\right\}. 
\end{equation}
The solution of this associativity equation, which is relatively well-known, is given as, \cite{Aczel66, Knuth10}:
\begin{equation}
	\label{BerC.10}
		\Theta\!\left[f\!\left(x, z\right)\right] = \Theta\!\left[f\!\left(x, y\right)\right] + \Theta\!\left[f\!\left(y, z\right)\right],
\end{equation}
where $\Theta\!\left(x\right)$ is some monotonic function.

If we substitute \eqref{BerC.8} into \eqref{BerC.10}, we get:
\begin{equation}
	\label{BerC.11}
		\Theta\!\left[h\!\left(\frac{x}{z}\right)\right] = \Theta\!\left[h\!\left(\frac{x}{y}\right)\right] + \Theta\!\left[h\!\left(\frac{y}{z}\right)\right]. 
\end{equation}
Letting $m\!\left(x\right) = \Theta\!\left[h\!\left(x\right)\right]$, we may rewrite \eqref{BerC.11} as:
\begin{equation}
	\label{BerC.11b}
		m\!\left(\frac{x}{z}\right) = m\!\left(\frac{x}{y}\right) + m\!\left(\frac{y}{z}\right). 
\end{equation}
Taking the derivative over $y$ \eqref{BerC.11b}, that is,
\begin{equation}
	\label{BerC.12}
		\frac{d}{dy} m\!\left(\frac{x}{z}\right) = \frac{d}{dy} \left[m\!\left(\frac{x}{y}\right) + m\!\left(\frac{y}{z}\right)\right], 
\end{equation} 
we obtain the ordinary differential equation:
\begin{equation}
	\label{BerC.13}
		0 = - \frac{x}{y^{2}} m'\!\left(\frac{x}{y}\right) + \frac{1}{z} m'\!\left(\frac{y}{z}\right), 
\end{equation} 
or, equivalently, 
\begin{equation}
	\label{BerC.14}
		\frac{x}{y} m'\!\left(\frac{x}{y}\right)  =  \frac{y}{z} m'\!\left(\frac{y}{z}\right). 
\end{equation} 
If we let $u = y/z$ and $c = x/z$, then we may rewrite \eqref{BerC.14} as
\begin{equation}
	\label{BerC.15}
		\frac{c}{u} m'\!\left(\frac{c}{u}\right)  =  u m'\!\left(u\right).
\end{equation} 
If we take the derivative of $c$ over \eqref{BerC.15}, that is, 
\begin{equation}
	\label{BerC.16}
		\frac{d}{dc} \frac{c}{u} m'\!\left(\frac{c}{u}\right)  =  \frac{d}{dc} u m'\!\left(u\right), 
\end{equation} 
then we obtain the ordinary differential equation:
\begin{equation}
	\label{BerC.17}
		\frac{1}{u} m'\!\left(\frac{c}{u}\right)  +  \frac{c}{u^{2}} m''\!\left(\frac{c}{u}\right) = 0. 
\end{equation} 
If we let $v = c/u$, then we may rewrite \eqref{BerC.17} as
\begin{equation}
	\label{BerC.18}
		m'\!\left(v\right)  +  v m''\!\left(v\right) = 0. 
\end{equation} 
The general solution of $m$ is given as:
\begin{equation}
	\label{BerC.19}
		m\!\left(v\right)  =  q \log v + C . 
\end{equation} 
where $q$ and $C$ are arbitrary constants.

It then follows from \eqref{BerC.8} and \eqref{BerC.19} that the function $f$, which adheres to the invariance desideratum \eqref{BerC.1} and the consistency desideratum \eqref{BerC.2}, has as its general solution:
\begin{equation}
	\label{BerC.20}
		f\!\left(x, y\right)  =  q \log \frac{x}{y} + C . 
\end{equation} 
where $q$ and $C$ are arbitrary constants.

From the boundary conditions \eqref{BerC.3} and \eqref{BerC.4}, it then follows, respectively, that $C = 0$ and that $q \geq 0$, which leaves us with the specific solution:
\begin{equation}
	\label{BerC.21}
		f\!\left(x, y\right)  =  q \log \frac{x}{y},  \quad \text{for } q \geq 0,
\end{equation} 
which is the Bernoulli utility function, originally proposed by Bernoulli in 1738, \cite{Bernoulli38}.

It follows that the Bernoulli utility function is the only function that adheres to the desiderata of unit invariance and consistency, respectively, \eqref{BerC.1} and \eqref{BerC.2}, and the boundary conditions that a zero change should lead to a zero perceived loss and that a perceived loss should be assigned a negative value, respectively, \eqref{BerC.3} and \eqref{BerC.4}. Any other utility function will be in violation with these fundamental desiderata and boundary conditions. 

The amazing thing here, at least to these authors' mind, is that so much was gained for so little. As one need not mention monetary gains and losses or changes in offered sensory stimuli in order to derive the Bernoulli utility function, or, equivalently, the Weber-Fechner law. 

We offer here the following observation, the functional equation \eqref{BerC.1} is, in its alternative interpretation, an invariance argument for two persons holding different amounts of initial wealth. Bernoulli, instead, uses variance arguments to arrive at his derivation. But both the variance and invariance arguments, as in the variance argument the wealth of person $A$ approaches the wealth of person $B$\footnote{See Appendix~\ref{BernoulliA}.}.

\section{The Criterion Of Choice}
\label{Choice}
In this section we will discuss the maximization of the sum of the lower and upper bounds, the third step in the Bayesian decision algorithm, as a criterion of choice. This criterion is highly non-intuitive, as it admits no interpretation. This is because it is a `corollary', or, to be more precise, an algebraic reshuffling, of a criterion which is intuitive, but less succinct in its expression.  

Let $D_{1}$ and $D_{2}$ be two decisions we have to choose from. Let $O_{i}$, for $i = 1, \dots, n$, and $O_{j}$, for $j = 1, \dots, m$, be the monetary outcomes associated with, respectively, decisions $D_{1}$ and $D_{2}$. 

In the Bayesian decision analysis, we first construct the two outcome distributions that correspond with these decisions:
\begin{equation}
	\label{eq.why.1}
	P\!\left(\left.O_{i}\right|D_{1}\right), \qquad  P\!\left(\left.O_{j}\right|D_{2}\right),
\end{equation}
where, if $n = m$, the outcomes $O_{i}$ and $O_{j}$ may or not may be equal for $i = j$.

We then proceed, by way of the Bernoulli utility function, or, equivalently, the Weber-Fechner law, to map utilities to the monetary outcomes $O_{i}$ and $O_{j}$ in \eqref{eq.why.1}. This leaves us with the utility probability distributions:
\begin{equation}
	\label{eq.why.2}
	P\!\left(\left.u_{i}\right|D_{1}\right), \qquad  P\!\left(\left.u_{j}\right|D_{2}\right).
\end{equation}

Our most primitive intuition regarding the utility probability distributions \eqref{eq.why.2} is that the decision which corresponds with the utility probability distribution which lies more to the right will also be the decision that promises to be the most advantageous. So, when making a decision we compare the positions of the utility probability distribution on the utility axis. This utility axis goes from minus infinity to plus infinity. Hence, the more-to-the-right criterion of choice.

Now, the confidence bounds of \eqref{eq.why.2}, say:
\begin{equation}
	\label{eq.why.3}
	\left[L\!B\!\left(\left.u_{i}\right|D_{1}\right), U\!B\!\left(\left.u_{i}\right|D_{1}\right)\right], \qquad  \left[L\!B\!\left(\left.u_{j}\right|D_{2}\right), U\!B\!\left(\left.u_{j}\right|D_{2}\right)\right],
\end{equation}
may provide us with a numerical handle on the concept of more-to-the-right.

For example, if we have that both
\begin{equation}
	\label{eq.why.4}
	L\!B\!\left(\left.u_{i}\right|D_{1}\right) > L\!B\!\left(\left.u_{j}\right|D_{2}\right),  \qquad U\!B\!\left(\left.u_{i}\right|D_{1}\right) > U\!B\!\left(\left.u_{j}\right|D_{2}\right).
\end{equation}
Then we will have an unambiguous preference for decision $D_{1}$ over decision $D_{2}$; seeing that under both the still probable worst and best case we will be better if we opt for $D_{1}$.

Likewise, if we have that either
\begin{equation}
	\label{eq.why.5}
	L\!B\!\left(\left.u_{i}\right|D_{1}\right) = L\!B\!\left(\left.u_{j}\right|D_{2}\right),  \qquad U\!B\!\left(\left.u_{i}\right|D_{1}\right) > U\!B\!\left(\left.u_{j}\right|D_{2}\right),
\end{equation}
or
\begin{equation}
	\label{eq.why.6}
	L\!B\!\left(\left.u_{i}\right|D_{1}\right) > L\!B\!\left(\left.u_{j}\right|D_{2}\right),  \qquad U\!B\!\left(\left.u_{i}\right|D_{1}\right) = U\!B\!\left(\left.u_{j}\right|D_{2}\right).
\end{equation}
Then, again, we will have an unambiguous preference for decision $D_{1}$ over decision $D_{2}$. In the constellation \eqref{eq.why.5}, we stand, all other things being equal, to be better of under the still probable best case scenario; while in the constellation \eqref{eq.why.6}, we stand, all other things being equal, to be less worse of under the still probable worst case scenario.

However, things become more ambiguous when, say, under decision $D_{1}$, we have to make a trade-off between either a gain in the upper bound and a loss in the lower bound
\begin{equation}
	\label{eq.why.7}
	L\!B\!\left(\left.u_{i}\right|D_{1}\right) < L\!B\!\left(\left.u_{j}\right|D_{2}\right), \qquad U\!B\!\left(\left.u_{i}\right|D_{1}\right) > U\!B\!\left(\left.u_{j}\right|D_{2}\right),
\end{equation}
or a gain in the lower bound and a loss in the upper bound
\begin{equation}
	\label{eq.why.8}
	L\!B\!\left(\left.u_{i}\right|D_{1}\right) > L\!B\!\left(\left.u_{j}\right|D_{2}\right), \qquad U\!B\!\left(\left.u_{i}\right|D_{1}\right) < U\!B\!\left(\left.u_{j}\right|D_{2}\right).
\end{equation}

We postulate here that a rational criterion of choice in the respective trade-off situations \eqref{eq.why.7} and \eqref{eq.why.8}, would be to pick that decision whose gain in either the lower or upper bound exceeds the loss in the corresponding upper or lower bound. 

So, if, say, under $D_{1}$ we stand to gain more in the still probable best case scenario than we stand to lose under the still probable worst case scenario, that is, \eqref{eq.why.7}:
\begin{equation}
	\label{eq.why.9}
	  L\!B\!\left(\left.u_{j}\right|D_{2}\right) - L\!B\!\left(\left.u_{i}\right|D_{1}\right) <  U\!B\!\left(\left.u_{i}\right|D_{1}\right) - U\!B\!\left(\left.u_{j}\right|D_{2}\right),
\end{equation}
then we will choose $D_{1}$ over $D_{2}$. Likewise, if under $D_{1}$ we stand to gain more in the still probable worst case scenario than we stand to lose under the still probable best case scenario, that is, \eqref{eq.why.8}:
\begin{equation}
	\label{eq.why.10}
	  L\!B\!\left(\left.u_{i}\right|D_{1}\right) - L\!B\!\left(\left.u_{j}\right|D_{2}\right) >  U\!B\!\left(\left.u_{j}\right|D_{2}\right) - U\!B\!\left(\left.u_{i}\right|D_{1}\right),
\end{equation}
then again we will choose $D_{1}$ over $D_{2}$.

Note that the gains and losses in this discussion pertain to gains and losses on the utility dimension, not on the monetary outcome dimension. On the utility dimension the phenomenon of loss aversion, that is, the phenomenon that monetary losses may weigh heavier than equal monetary qains, has already been accounted for. Stated differently, the utility scale is a linear loss-aversion corrected scale for the moral value of monies.

Now, if we look at the scenarios \eqref{eq.why.7} and \eqref{eq.why.8}, and the corresponding postulated rational, because intuitive, criteria of choice \eqref{eq.why.9} and \eqref{eq.why.10}, then we see that we will choose $D_{1}$ over $D_{2}$ whenever we have that
\begin{equation}
	\label{eq.why.11}
	L\!B\!\left(\left.u_{i}\right|D_{1}\right) + U\!B\!\left(\left.u_{i}\right|D_{1}\right) > L\!B\!\left(\left.u_{j}\right|D_{2}\right) + U\!B\!\left(\left.u_{j}\right|D_{2}\right).
\end{equation}
Moreover, this single criterion of choice is also consistent with the choosing of $D_{1}$ over $D_{2}$ in the scenarios \eqref{eq.why.4}, \eqref{eq.why.5}, and \eqref{eq.why.6}. 

This, then, is the rationale behind the non-intuitive, because it admits no interpretation, criterion of choice that we should maximize the sum of the lower and upper bounds of the utility probability distributions, in order to come to the optimal decision\footnote{Instead of the criterion \eqref{eq.why.11}, one may also use a lower bound maximization, that is, choose $D_{1}$ whenever
\[
	L\!B\!\left(\left.u_{i}\right|D_{1}\right) > L\!B\!\left(\left.u_{j}\right|D_{2}\right),
\]
A possible lower bound maximizer might be the regulator who is not that interested in a bank's potential profit, but only has an eye for the potential catastrophic losses that, were they to materialize, could destabilize the entire financial system. Or,alternatively, one may use an upper bound maximization, that is, choose $D_{1}$ whenever
\[
	U\!B\!\left(\left.u_{i}\right|D_{1}\right) > U\!B\!\left(\left.u_{j}\right|D_{2}\right),
\]
where the upper bound maximizer is the banker that tries to maximize his yearly bonus.}.

Note that if the decision inequality \eqref{eq.why.11} goes to an equality:
\begin{equation}
	\label{eq.why.12}
	L\!B\!\left(\left.u_{i}\right|D_{1}\right) + U\!B\!\left(\left.u_{i}\right|D_{1}\right) = L\!B\!\left(\left.u_{j}\right|D_{2}\right) + U\!B\!\left(\left.u_{j}\right|D_{2}\right).
\end{equation}
Then we have that we will be undecided when it comes to the decisions $D_{1}$ and $D_{2}$. 

Also note that for $k$-sigma bounds \eqref{eq.why.3} translates to
\begin{equation}
	\label{eq.why.13}
	E\!\left(\left.u_{i}\right|D_{1}\right) \pm k \: \text{std}\!\left(\left.u_{i}\right|D_{1}\right), \qquad  E\!\left(\left.u_{j}\right|D_{2}\right) \pm k \: \text{std}\!\left(\left.u_{j}\right|D_{2}\right), 
\end{equation}
which, if substituted in \eqref{eq.why.11}, gives the inequality
\begin{equation}
	\label{eq.why.15}
	2 E\!\left(\left.u_{i}\right|D_{1}\right) > 2 E\!\left(\left.u_{j}\right|D_{2}\right),
\end{equation}
which brings us right back to Bernoulli's expected utility theory, as proposed in 1738, \cite{Bernoulli38}, in which it is stated that the expectation value of the utility probability distribution should be maximized. 

Nonetheless, the criterion of choice, that the sum of the upper and lower bound should be maximized, as proposed here, will deviate from Bernoulli's initial 1738 proposal when the $k$-sigma intervals overshoots either its minimal or maximal value of the utility probability distribution.

Let $a$ and $b$, respectively, be the minimal and maximal values of a given utility probability distribution. Then we may identify two additional symmetry breaking cases, relative to \eqref{eq.why.15}:
\begin{equation}
	\label{eq.why.15b}
	L\! B\!\left(u\right) + U\! B\!\left(u\right) = 
	\begin{cases} 2 \: E\!\left(u\right), \qquad  &L\! B\!\left(u\right) > a, \quad U\! B\!\left(u\right) < b \\
	a + U\! B\!\left(u\right), \qquad &L\! B\!\left(u\right) < a, \quad U\! B\!\left(u\right) < b \\
	L\! B\!\left(u\right) + b, \qquad   &L\! B\!\left(u\right) > a, \quad U\! B\!\left(u\right) > b 
	\end{cases}
\end{equation}
These two last symmetry breaking cases, where the $k$-sigma intervals overshoot the minimal and maximal values of a given utility probability distribution are non-trivial, as they give rise to the specific $S$-shape of the fair probabilities in certainty bets\footnote{Certainty bets will be discussed in the section following the next one.}.

Another instance where we will deviate from Bernoulli's proposal is when we put an explicit premium on either caution or opportunity. If we take as the lower and upper bounds, whose sum is to be maximized:
\begin{equation}
	\label{eq.why.16}
	L\!B\!\left(u\right) = E\!\left(u\right) - a \: \text{std}\!\left(u\right),
\end{equation}	
and
\begin{equation}
	\label{eq.why.17}
	U\!B\!\left(u\right) = E\!\left(u\right) + b \: \text{std}\!\left(u\right).
\end{equation}
Then \eqref{eq.why.16} and \eqref{eq.why.17} sum to:
\begin{equation}
	\label{eq.why.18}
	L\!B\!\left(u\right) + U\!B\!\left(u\right) = 2 \: E\!\left(u\right) + \left(b - a\right) \: \text{std}\!\left(u\right).
\end{equation}
If in \eqref{eq.why.18} we let $a > b$, then we put a premium caution; if we set $b > a$, then we put a premium on opportunity; and if we let $a = b$, then we have an equal trade-off between caution and opportunity taking.

\section{The Ellsberg Paradox}
\label{chapter2.3}
In this section we will demonstrate how to construct a non-trivial outcome distribution, by way of the product and sum rules, and how to map outcomes to their corresponding utilities by way of the product and sum rules. 

Ellsberg found that the willingness to bet on an uncertain event depends not only on the degree of uncertainty but also on its source, \cite{Ellsberg61}. He observed that people prefer to bet on an urn containing equal numbers of red and green balls, rather than on an urn that contains red and green balls in unknown proportions. Ellsberg called this observed phenomenon source dependence.

Tversky and Kahneman \cite{Tversky92}, state that source dependence constitutes one of the minimal challenges that must be met by any adequate descriptive theory of choice. As our theory of choice is Bayesian, we will proceed to give a Bayesian treatment of this phenomenon\footnote{Note that we do not claim that decision makers will derive the Bayesian equations which will follow \textit{ad verbatim}. Rather, we state that Bayesian inference is common sense amplified, having a much higher probability resolution than our human brains can ever hope to achieve. So, being common sense amplified, the \textit{results} of the Bayesian analysis should be commensurate with our intuitions; rather than the analysis itself. See also Appendix~\ref{Fred}.}. 

\subsection{Constructing outcome probability distributions.}
Say, we have a large urn consisting of 1000 balls of which 500 are red and 500 green. We tell our subject that of the $N=1000$ balls $R=500$ are red and $N-R=500$ green, and that he is to draw a ball $n = 100$ times. After each draw he will get a dollar if the ball is red and nothing if the ball is green, after which the ball is to be put back in the urn. The subject is also told that for the privilege to partake in this bet an entrance fee of 50 dollars is to be paid. 

\begin{figure}[h]
		\centering
				\includegraphics[width=0.60\textwidth] {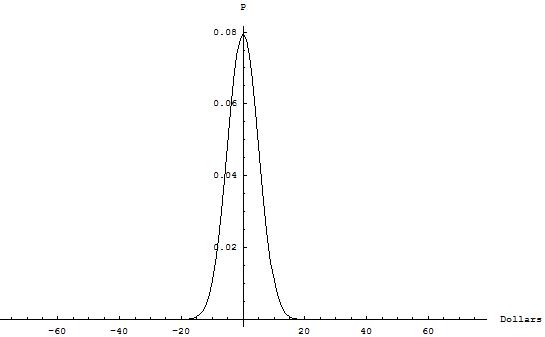} 
		\caption{Probability outcome distribution for bet 1}
		\label{fig:P1.3.1}
	\end{figure}

The probability of drawing $r$ red balls in the first bet $D_{1}$ may be modeled by way of a binomial distribution:
\begin{equation}
	\label{eq.P1.3.1}
	p\!\left(\left.r\right|n, R, N, D_{1}\right) = \frac{n!}{r! \left(n-r\right)!} \left(\frac{R}{N}\right)^{r} \left(1-\frac{R}{N}\right)^{n-r}.
\end{equation}
Now as the net return, say, $o$ is in dollars, having as its value number of red balls minus entrance fee, we have
\begin{equation}
	\label{eq.P1.3.2}
	o = r - 50,
\end{equation}
where, as there are $n = 100$ draws, $-50 \leq o \leq 50$. We then make a simple change of variable, using \eqref{eq.P1.3.2},
\begin{equation}
	\label{eq.P1.3.3}
	r = o + 50
\end{equation}
and substitute \eqref{eq.P1.3.3} into \eqref{eq.P1.3.1}, so as to get the probability function of the net return
\begin{equation}
	\label{eq.P1.3.4}
	p\!\left(\left.o\right|n, R, N, D_{1}\right) = \frac{n!}{\left(o+50\right)! \left(n-o-50\right)!} \left(\frac{R}{N}\right)^{o+50} \left(1-\frac{R}{N}\right)^{n-o-50}
\end{equation}
The probability distribution of the net return for bet $D_{1}$, that is, 
\[
	p\!\left(\left.o\right|n=100, R=500, N=1000, D_{1}\right),
\]
then can be plotted as, Figure~\ref{fig:P1.3.1}:\\

\noindent This probability distribution has a mean and a standard deviation of, respectively,
\begin{equation}
	\label{eq.P1.3.5}
	E\!\left(\left.o\right|D_{1}\right) = 0, \qquad \text{std}\!\left(\left.o\right|D_{1}\right) = 5.
\end{equation}

For the second bet, we tell our subject that the urn holds $N=1000$ balls which are either red or green. Again, for every red ball drawn there will be a dollar payout. There will be $n=100$ draws and after each draw the ball is to be placed in the urn again. The entrance fee of the bet is 50 dollars. However, we are also told that the number of red balls is neither zero nor thousand\footnote{We could let the range of the red balls be $0 \leq R \leq N$. But then the resulting outcome probability distribution would become unwieldy, because of this added structure.}, that is, $0 < R < N$, thus, precluding the certainty outcomes of $o = -50$ and $o = 50$ for, respectively, $R = 0$ and $R = 1000$. 

As the subject does not know the actual number of red balls $R$, the Bayesian thing to do is to weigh the probability of drawing $r$ red balls over all plausible values of $R$, the total number red balls in the urn. 

Based on the background information, we assign an uniform prior probability distribution to $R$, the unknown number of red balls in the urn:
\begin{equation}
	\label{eq.P1.3.6}
	p\!\left(R\right) = \frac{1}{N-1},
\end{equation}
where $R = 1, \ldots, N - 1$. So, by way of the product and the generalized sum rules, \eqref{eq.P1.2.3} and \eqref{eq.P1.2.4}, the probability of drawing $r$ red balls in the second bet $D_{2}$ translates to
\begin{equation}
	\label{eq.P1.3.7}
	p\!\left(\left.r\right|n, N, D_{2}\right) = \sum_{R = 1}^{N-1} p\!\left(\left.r, R\right|n, N, D_{2}\right) = \sum_{R = 1}^{N-1} p\!\left(R\right) p\!\left(\left.r\right|n, R, N, D_{2}\right).
\end{equation}
Again making a change of variable from the number of red balls drawn $r$ to the net return $o$, we substitute \eqref{eq.P1.3.1} and \eqref{eq.P1.3.6} into \eqref{eq.P1.3.7} and make a change of variable by way of \eqref{eq.P1.3.3}. This results in the probability distribution of the net return $o$:
\begin{equation}
	\label{eq.P1.3.8}
	p\!\left(\left.o\right|n, N, D_{2}\right) = \sum_{R = 1}^{N-1} \frac{1}{N-1} \frac{n!}{\left(o+50\right)!\left(n-o-50\right)!}  \left(\frac{R}{N}\right)^{o+50} \left(1-\frac{R}{N}\right)^{n-o-50}.
\end{equation}
The probability distribution of the net return for bet $D_{2}$, that is,, 
\[
	 p\!\left(\left.o\right|n=100, N=1000, D_{2}\right),
\]
then can be plotted as, Figure~\ref{fig:P1.3.2}:

\begin{figure}[h]
		\centering
			\includegraphics[width=0.60\textwidth]{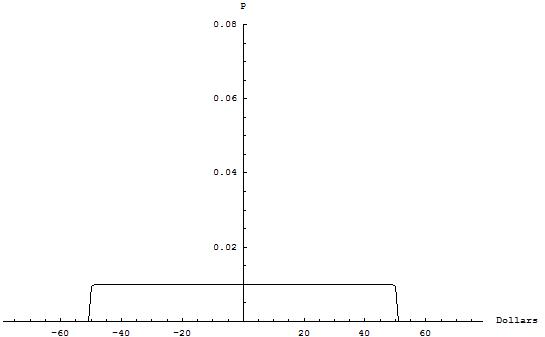}
		\caption{Probability outcome distribution for bet 2}
		\label{fig:P1.3.2}
	\end{figure}
	
\newpage	
	
\noindent This probability distribution has a mean and standard deviation of, respectively,
\begin{equation}
	\label{eq.P1.3.9}
	E\!\left(\left. o\right|D_{2}\right) = 0, \qquad  \text{std}\!\left(\left. o\right|D_{2}\right) = 29.
\end{equation}

\begin{figure}[h]
		\centering
			\includegraphics[width=0.60\textwidth]{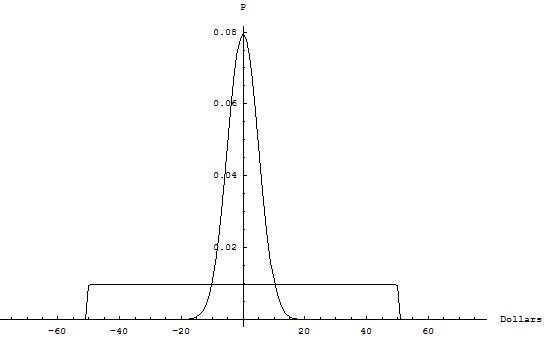}
		\caption{Probability outcome distributions for bets 1 and 2}
		\label{fig:P1.3.3}
	\end{figure}

In Figure~\ref{fig:P1.3.3}, we give both probability distributions, Figures 6.1 and 6.2, together. \\

\subsection{Constructing utility probability distributions.}
Let $D_{1}$ stand for the decision to choose the first Ellsberg bet and $D_{2}$ for the decision to choose the second Ellsberg bet. Then, by substituting, where appropriate, the values $n=100$, $R=500$, $N=1000$ into \eqref{eq.P1.3.4} and \eqref{eq.P1.3.8}, we may obtain the corresponding outcome probability distributions:
\begin{align}
	\label{eq.P1.4.3}
	&p\!\left(\left.o\right|n = 100, R = 500, N = 1000, D_{1}\right) = \frac{100!}{\left(o+50\right)! \left(50-o\right)!} \left(\frac{1}{2}\right)^{100}, \notag \\
	\\
		&p\!\left(\left.o\right|n = 100, N = 1000, D_{2}\right) =  \frac{100!}{\left(o+50\right)!\left(50-o\right)!} \sum_{R = 1}^{999}   \frac{1}{999} \left(\frac{R}{1000}\right)^{o+50} \left(1-\frac{R}{1000}\right)^{50-o}. \notag 
\end{align}

In order to map the outcomes $o$ in \eqref{eq.P1.4.3} to values on the utility dimension $u$, we introduce the conditional probability distribution $p\!\left(\left.u\right|o\right)$. Combining \eqref{eq.P1.4.3} with $p\!\left(\left.u\right|o\right)$, by way of the product rule \eqref{eq.P1.2.3}, and marginalizing over the possible outcomes $o$, by way of the generalized sum rule \eqref{eq.P1.2.4}, we may get the utility probability distribution of interest, that is,
\begin{equation}
	\label{eq.P1.5.1b}
	p\!\left(\left.u\right|D_{i}\right) = \sum_{o} p\!\left(\left. u, o\right|D_{i}\right) = \sum_{o} p\!\left(\left.u\right|o\right) p\!\left(\left.o\right|D_{i}\right).
\end{equation}

Now, if we assign as our utility function the Bernoulli utility function, then the conditional utility probability distribution to be employed in the Ellsberg example is:
\begin{equation}
	\label{eq.P1.4.4}
p\!\left(\left.u\right|o, m, q\right) =  \begin{cases}
										1,	&u = q \log \frac{m + o}{m}   \\
										0,	&u\neq q \log \frac{m + o}{m}
										\end{cases}
\end{equation}
for $o = -50, -49, \ldots, 50$, and where $m$ is our initial wealth and $q$ is the scaling constant of the Bernoulli utility function. This probability distribution takes us from the $o$-dimension, which is the dimension of the monetary outcomes, to the $u$-dimension, which is the dimension of the moral value of these monetary outcomes. 

From \eqref{eq.P1.4.4}, we see that every outcome $o$ admits only one utility value $u$; that is, \eqref{eq.P1.4.4} is of the Dirac delta form: 
\begin{equation}
	\label{eq.P1.6.2q}
	p\!\left(\left.u\right|o, q, m\right) = \delta\!\left(u - q \log \frac{m + o}{m}\right),
\end{equation}
where $\delta$ is the Dirac delta function for which
\begin{equation}
	\label{eq.P1.6.2c}
	\delta\!\left(u-c\right) du = 
										\begin{cases}
										1,	&u = c   \\
										0,	&u \neq c
										\end{cases}
\end{equation}
Because of \eqref{eq.P1.6.2c}, we have that
\begin{equation}
	\label{eq.P1.6.2c2}
	\int \delta\!\left(u-c\right) f\!\left(u\right) du = f\!\left(c\right).
\end{equation}

Applying both  \eqref{eq.P1.5.1b} and the Dirac delta \eqref{eq.P1.6.2q} to the outcome distributions \eqref{eq.P1.4.3}, we obtain the utility probability distributions:
\begin{align}
	\label{eq.P1.4.5}
	p\!\left(\left.u\right|m, q, D_{1}\right) &= \sum_{o = -50}^{50} \delta\!\left(u - q \log \frac{m + o}{m}\right) p\!\left(\left.o\right| D_{1}\right), \notag \\
	\\
	p\!\left(\left.u\right|m, q, D_{2}\right) &= \sum_{o = -50}^{50} \delta\!\left(u - q \log \frac{m + o}{m}\right)  p\!\left(\left.o\right| D_{2}\right),\notag 
\end{align}
where in the probability distributions we only conditionalize on that which is not yet specified. 

The $k$th-order moments of the utility probability distributions \eqref{eq.P1.4.5} may be evaluated by way of the integrals: 
\begin{equation}
	\label{eq.P1.4.6}
	E\!\left(\left.u^{k}\right|D_{i}\right) = \int u^{k} \; p\!\left(\left.u\right|m, q, D_{i}\right) du.
\end{equation}

\subsection{Applying the criterion of choice.}
If we assume a modest intial wealth, that is, monthly expendable income, of two-hundred dollars, that is, $m = 200$. Then the mean, standard deviation, and skewness of the utility probability distributions \eqref{eq.P1.4.5} are given as\footnote{The first two cumulants of the utility probability distributions \eqref{eq.P1.4.5} may be computed by way of \eqref{eq.P1.4.6} and the identities, \cite{Hall92}:
\begin{equation}
	\mu = E\!\left(\left.u\right|D_{i}\right), \qquad \sigma = \sqrt{E\!\left(\left.u^{2}\right|D_{i}\right) - \left[E\!\left(\left.u\right|D_{i}\right)\right]^{2}}. \notag
\end{equation}}
\begin{equation}
	\label{eq.ells.1}
	E\!\left(\left.u\right|q, D_{1}\right) = -0.0003 \; q, \qquad \text{std}\!\left(\left.u\right|q, D_{1}\right) = 0.025 \: q, 
\end{equation}
and  
\begin{equation}
	\label{eq.ells.2}
	E\!\left(\left.u\right|q, D_{2}\right) = -0.0108 \: q, \qquad \text{std}\!\left(\left.u\right|q, D_{2}\right) = 0.1479 \: q. 
\end{equation}
The minimum and maximum values, respectively, $a$ and $b$, of the utility probability distributions are given as, \eqref{eq.P1.5.4b} and problem statement:
\begin{equation}
	\label{eq.ells.2b}
	a = q \log \frac{200 - 50}{200} = - 0.2877 q, \qquad b = q \log \frac{200 - 50 + 100}{200} = 0.2231 q.
\end{equation}

In the following, we will use 1-sigma intervals. Seeing that for these confidence intervals no overshoot occurs, we use for the sum of the lower and upper bounds of the utility probability distributions the first identity of \eqref{eq.why.15b}. This gives for \eqref{eq.ells.1}:
\begin{equation}
	\label{eq.ells.3}
	L\!B\!\left(\left.u\right|q, D_{1}\right) + U\!B\!\left(\left.u\right|q, D_{1}\right) = 2 E\!\left(\left.u\right|q, D_{1}\right) = -0.0006 \: q,
\end{equation}
and for \eqref{eq.ells.2}:
\begin{equation}
	\label{eq.ells.4}
	L\!B\!\left(\left.u\right|q, D_{2}\right) + U\!B\!\left(\left.u\right|q, D_{2}\right) = 2 E\!\left(\left.u\right|q, D_{1}\right) = -0.0216 \: q.
\end{equation}

If we compare \eqref{eq.ells.3} and \eqref{eq.ells.4}, then we find that decision $D_{1}$ is more advantageous than decision $D_{2}$:
\begin{equation}
	\label{eq.ells.5}
	L\!B\!\left(\left.u\right|q, D_{1}\right) + U\!B\!\left(\left.u\right|q, D_{1}\right) > L\!B\!\left(\left.u\right|q, D_{2}\right) + U\!B\!\left(\left.u\right|q, D_{2}\right),
\end{equation}	
since
\begin{equation}
	\label{eq.ells.6}
	-0.0006 \; q > -0.0216 \; q,
\end{equation}	
for any positive scaling constant $q$. 

As an aside the unknown scaling constant $q$ will always fall away in the decision theoretical inequalities. This is because the mean and standard deviations of the utility probability distributions are linear in $q$. Stated differently, it may be checked that, for the stochastics $X$ and $Y$, and a positive constant $q$, \cite{Lindgren93}: 
\begin{equation}
	\label{eq.P1.D.11}
	E\!\left(q X\right) = q \: E\!\left(X\right), \qquad \text{std}\!\left(q X\right) = q \; \text{std}\!\left(X\right).
\end{equation}
So, unless we want to interpret our decisions in terms of net utility gain, as we do below, then we may, without any loss of generality, set the unknown scaling constant $q$ to one. 

Now, if we are willing to commit ourselves to the scaling constant value of $q = 100$ as the Weber constant for monetary stimuli, \eqref{eq.P1.10.2b}. Then we may, by way of \eqref{eq.why.10}, interpret \eqref{eq.ells.5} as follows. 

Under decision $D_{1}$, betting on an urn with an equal number of red and green balls, there is a gain of 13.34 utiles, in terms of loss mitigation, and a loss of 11.24 utiles, in terms of gain reduction, relative to decision $D_{2}$, betting on an urn with an unknown proportion of red balls. This makes $D_{1}$, with a net utility gain of 2.10 utiles, more attractive a choice than $D_{2}$. 

\begin{figure}[h]
		\centering
			\includegraphics[width=0.60\textwidth]{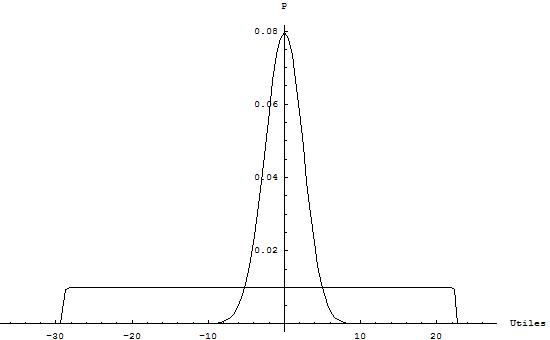}
		\caption{Probability utility distributions for bets 1 and 2}
		\label{fig:P1.3.4}
	\end{figure}

The setting of the scaling constant $q$ also allows us to plot the utility probability distributions under $D_{1}$ and $D_{2}$, Figure~\ref{fig:P1.3.4}.

We summarize, in a Bayesian decision theoretical analysis we first construct the outcome probability distributions on a monetary unit scale under the decisions $D_{1}$ and $D_{2}$, Figure~\ref{fig:P1.3.3}. We then construct, by way of the Bernoulli utility function, the corresponding utility probability distributions on an (un)scaled utile scale, Figure~\ref{fig:P1.3.4}. We then compare, for a given decision, the gain/loss in the lower bound relative to the corresponding loss/gain in the upper bound; or, equivalently, in the algebraic sense of the word, we compare the sums of the upper and lower bounds under the decisions $D_{1}$ and $D_{2}$. 

The results of this Bayesian decision theoretical analysis is in correspondence with the Ellsberg finding that people prefer to bet on an urn containing equal numbers of red and green balls, rather than on an urn that contains red and green balls in unknown proportions, \cite{Ellsberg61}.

In closing, we may envisage decision problems in which we are uncertain regarding the actual utility of a given outcome $o$. Such an occasion may arise when we do not know the initial wealth $m$ of the subject under investigation. In those cases we will want to assign probability distributions $p\!\left(\left.u\right|o\right)$ less dogmatic than the Dirac delta \eqref{eq.P1.6.2q} to our utilities.

\section{The Psychological Certainty Effect, Part I}
\label{KT1}
Certainty bets admit the following structure. Let $O_{c}$ and $O_{u}$, respectively, be the certainty and the uncertainty outcomes, where $O_{c} < O_{u}$. The uncertainty outcome $O_{u}$ has a probability of $p$ of being realized. In the case that $O_{u}$ is not realized the outcome will be zero, and the probability corresponding with this outcome is $1-p$. The certainty outcome $O_{c}$ is certain and, hence, has a probability of one.

So, the outcome probability distributions for certainty bets are 
\begin{equation}
	\label{eq.P1.10.1q}
	 p\!\left(\left.O_{i}\right|D_{1}\right) = \begin{cases} 
	 										   p, \qquad &O_{1} = O_{u} \\
										     1-p, \qquad &O_{2} = 0
											\end{cases}
\end{equation} 
and
\begin{equation}
	\label{eq.P1.10.2q}
	 p\!\left(\left.O_{j}\right|D_{2}\right) = \begin{cases} 
	 											1.0, \qquad O_{1} = O_{c} 
											\end{cases}
\end{equation}

In what follows we will assume, initially, for simplicity's sake, a linear utility for monetary outcomes. This assumption corresponds with an initial wealth $m$ that vastly exceeds any increment $\Delta m$. This can be seen as follows. The logarithmic function admits the series expansion:
\begin{equation}
	\label{eq.P1.10.2qq}
	 \log\left(1 + x\right) = x -\frac{x^{2}}{2} + \frac{x^{3}}{3} -\frac{x^{4}}{4} + \cdots.
\end{equation} 
So, as $m /\Delta m$ tends to zero, we have that the Bernoulli utility function, or, equivalently, the Weber-Fechner law, \eqref{eq.P1.5.4b} tends to, \eqref{eq.P1.10.2qq}:
\begin{equation}
	\label{eq.P1.10.2c}
	 q \log \frac{m + \Delta m}{m} = q \log\left(1 + \frac{\Delta m}{m}\right) \rightarrow q \frac{\Delta m}{m},
\end{equation} 
as can be seen in Figure~\ref{fig:P1.5.2}.

Fairness is defined as the decision theoretical equality:
\begin{equation}
	\label{eq.P1.10.30}
	 LB\!\left(\left.u\right| p, D_{1}\right) + UB\!\left(\left.u\right| p, D_{1}\right) =  LB\!\left(\left.u\right| D_{2}\right) + UB\!\left(\left.u\right| D_{2}\right),
\end{equation}
where $D_{1}$ and $D_{2}$ correspond with the choosing of, respectively, the uncertainty and certainty bets. Fairness finds its expression in the equality \eqref{eq.P1.10.30} because this equality guarantees that any loss/gain in the lower bound will be offset by a commensurate gain/loss in the upper bound\footnote{See Section~\ref{Choice}.}, for both $D_{1}$ and $D_{2}$.

If, for a certainty bet having positive outcomes, we solve \eqref{eq.P1.10.30} for the fair probability $p$, assuming a linear utility for money, and taking care to take into account any symmetry breaking conditions that may occur, \eqref{eq.why.15b}, we find that the fair probability $p$ maps to the outcome intervals   
\begin{equation}
	\label{eq.P1.10.31}
	\left(0, 2 O_{c}\right),		\qquad \text{for } O_{c} \leq \frac{O_{u}}{2},
\end{equation}
which is intuitively fair for both the takers and the providers of decision $D_{1}$, relative to the certainty offer of $O_{c}$, and
\begin{equation}
	\label{eq.P1.10.32}
	\left(2 O_{c} - O_{u}, O_{u}\right),		\qquad \text{for } O_{c} > \frac{O_{u}}{2},
\end{equation}
which is intuitively fair for both the takers and providers of decision $D_{1}$, relative to the certainty offer of $O_{c}$.

If for an uncertainty pay out of either $0$ or $O_{u} = 5000$, we plot the solution of \eqref{eq.P1.10.30} for the fairness probability $p$, assuming a linear utility for monetary outcomes, as a function of the certainty outcome $O_{c}$, we obtain Figure~\ref{fig:risk_seeking_1e}:

\begin{figure}[h]
	\centering
		\includegraphics[width=0.60\textwidth]{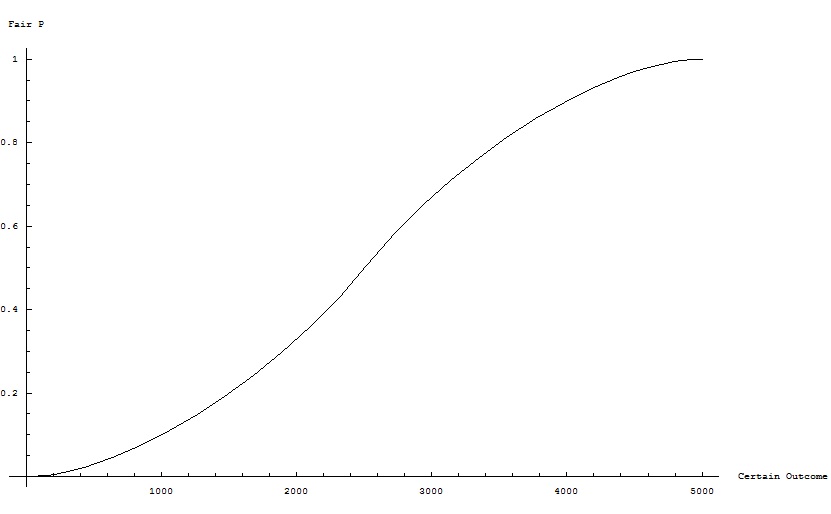}
	\caption{Fair probability as function of certain outcome, positive outcomes}
	\label{fig:risk_seeking_1e}
\end{figure}

If we again solve \eqref{eq.P1.10.30} for the fairness probability $p$, as a function of the certainty outcome $O_{c}$, assuming a linear utility for monetary outcomes, but now neglecting any symmetry breaking conditions, and add this expected utility theory solution to Figure~\ref{fig:risk_seeking_1e}, we obtain Figure~\ref{fig:risk_seeking_1f}:

\begin{figure}[h]
	\centering
		\includegraphics[width=0.60\textwidth]{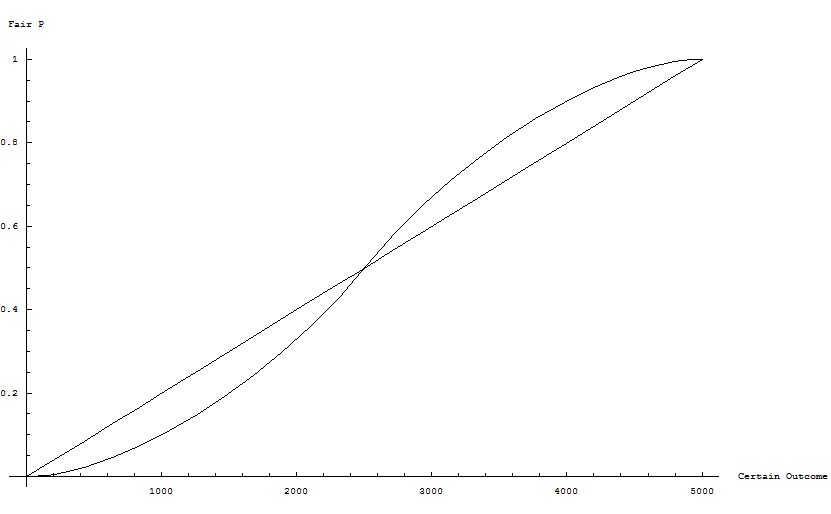}
		\caption{Fair probability for sum of bounds and expectation value maximization, positive outcomes}
	\label{fig:risk_seeking_1f}
\end{figure}

\newpage

If, for a certainty bet having negative outcomes, we solve \eqref{eq.P1.10.30} for the fair probability $p$, assuming a linear utility for money, and taking care to take into account any symmetry breaking conditions that may occur, \eqref{eq.why.15b}, we find that the fair probability $p$ maps to the outcome intervals   
\begin{equation}
	\label{eq.P1.10.31b}
	\left(2 O_{c}, 0\right),		\qquad \text{for } O_{c} \geq \frac{O_{u}}{2},
\end{equation}
which is intuitively fair for both the takers and the providers of decision $D_{1}$, relative to the certainty offer of $O_{c}$, and
\begin{equation}
	\label{eq.P1.10.32b}
	\left(O_{u}, 2 O_{c} - O_{u}\right),		\qquad \text{for } O_{c} < \frac{O_{u}}{2},
\end{equation}
which is intuitively fair for both the takers and providers of decision $D_{1}$, relative to the certainty offer of $O_{c}$.

If for an uncertainty pay out of either $0$ or $O_{u} = -5000$, we plot the solution of \eqref{eq.P1.10.30} for the fairness probability $p$, assuming a linear utility for monetary outcomes, as a function of the negative certainty outcome $O_{c}$, we obtain Figure~\ref{fig:risk_seeking_1h}:

\begin{figure}[h]
	\centering
		\includegraphics[width=0.60\textwidth]{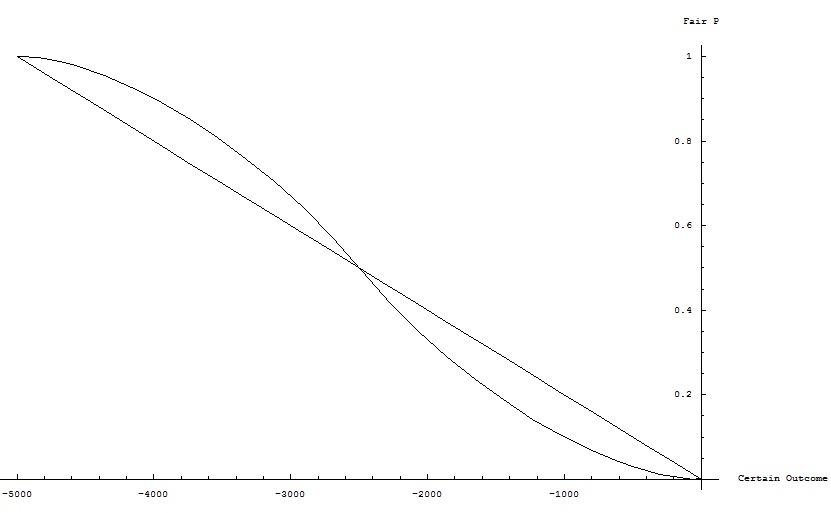}
		\caption{Fair probability for sum of bounds and expectation value maximization, negative outcomes}
	\label{fig:risk_seeking_1h}
\end{figure}

If we rescale the $x$-axes of Figures~\ref{fig:risk_seeking_1f} and~\ref{fig:risk_seeking_1h} as the ratio $O_{c}/O_{u}$, where $\left|O_{c}\right| \leq \left|O_{u}\right|$, and reverse the axes, we obtain the alternative Figure~\ref{fig:risk_seeking_1g}:

\begin{figure}[h]
	\centering
		\includegraphics[width=0.60\textwidth]{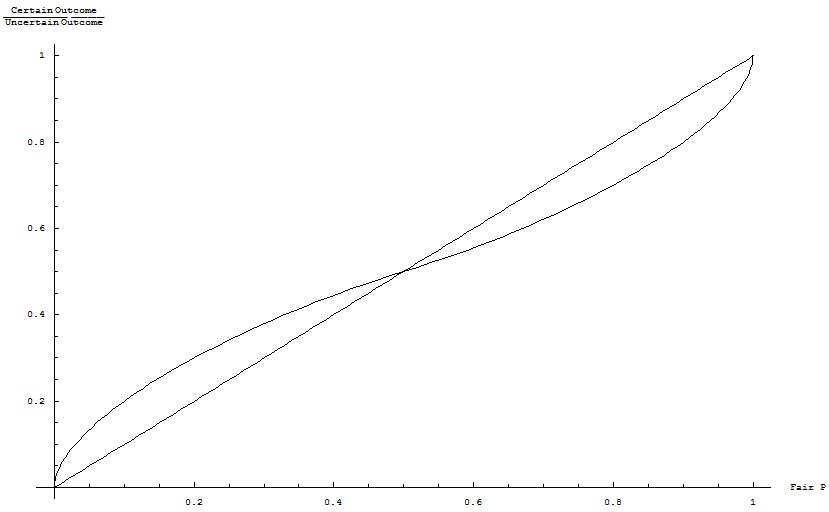}
		\caption{Rescaled and rotated figure for positive and negative values}
	\label{fig:risk_seeking_1g}
\end{figure}

Now, those readers who are familiar with the cumulative prospect theory may recognize in Figure~\ref{fig:risk_seeking_1g}, Kahneman and Tversky's Figures~1, 2, and~3 of their \cite{Tversky92}. 

But Kahneman and Tversky obtained their figures not from first principles, as is done here, but through experimentation, in which subjects where asked to decide on certainty bets of the type we discussed in the previous section. So, it would seem that Kahneman and Tversky, inadvertently, for they are outspoken anti-Bayesian\footnote{See Appendix~\ref{Represent}.}, have provided the Bayesian decision theory with a very strong supporting contact. 

Kahneman and Tversky see in the empirical observation of the typical $S$-curve of Figure~\ref{fig:risk_seeking_1g} another justification\footnote{See Appendix~\ref{Non-Linear Preferences}, for a discussion of the initial justification.} for their probability weighing functions\footnote{Note that Kahneman and Tversky's $\gamma$ and $\delta$ are not our $\gamma$ and $\delta$.},
\begin{equation}
	\label{eq.P1.10.33a}
	w^{+}\!\left(p\right) = \frac{p^{\gamma}}{p^{\gamma} + \left(1- p\right)^{\frac{1}{\gamma}}},
\end{equation}       
and
\begin{equation}
	\label{eq.P1.10.33b}
	w^{-}\!\left(p\right) = \frac{p^{\delta}}{p^{\delta} + \left(1- p\right)^{\frac{1}{\delta}}},
\end{equation}     
which over weighs small probabilities and under weighs large probabilities. Moreover, Kahneman and Tversky offer up the implied under weighing of small probabilities, in order to explain the general popularity of lotteries and insurances. 

We, on the other hand, see in the empirical observation of the typical $S$-curve of Figure~\ref{fig:risk_seeking_1g} a confirmation of the non-triviality of the proposed criterion of choice \eqref{eq.why.15b}. 

If we drop the assumption of a linear utility of monetary outcomes in the neighborhood of $-5000 < \Delta m < 5000$, and for initial wealths of $m = 1000$ and $m = 6000$ for certainty bets involving, respectively, positive and negative outcomes. Then we may assign, by way of the Bernoulli law, \eqref{eq.P1.5.4b}, utilities to the monetary outcomes. By doing so, we obtain the following fairness ratio outcomes for a given probability $p$ of the uncertain proposition, Figures~\ref{fig:risk_seeking_1i} and~\ref{fig:risk_seeking_1j}:  

\begin{figure}[h]
	\centering
		\includegraphics[width=0.60\textwidth]{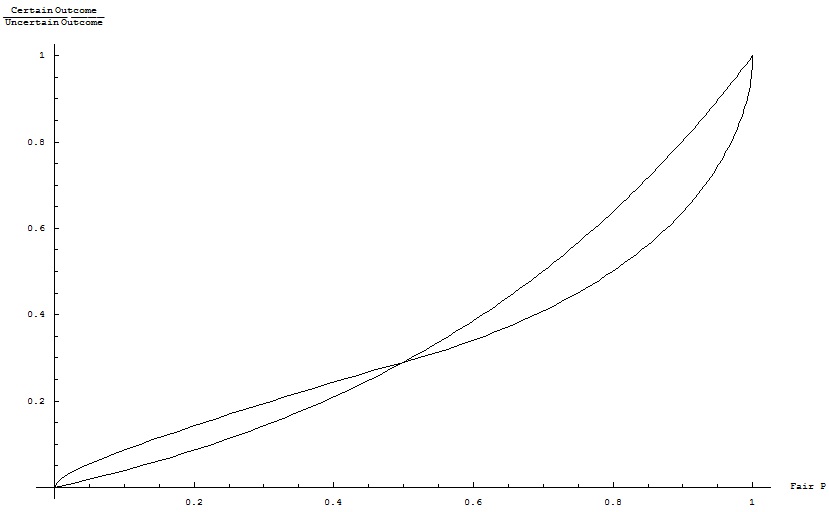}
	\caption{Rescaled and rotated figure for positive outcomes}
	\label{fig:risk_seeking_1i}
\end{figure}

and

\begin{figure}[h]
	\centering
		\includegraphics[width=0.60\textwidth]{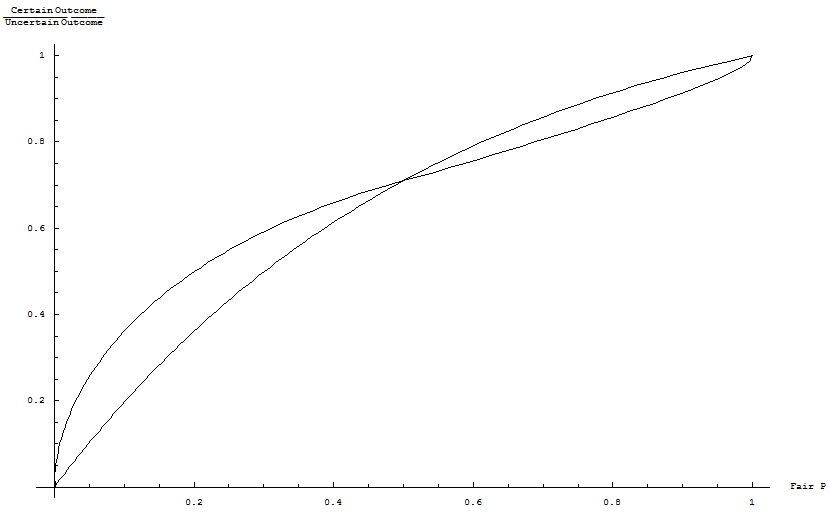}
		\caption{Rescaled and rotated figure for negative outcomes}
	\label{fig:risk_seeking_1j}
\end{figure}

Comparing Figures~\ref{fig:risk_seeking_1i} and~\ref{fig:risk_seeking_1g}, we see that by taking into account the initial wealth $m$, through the Bernoulli law, \eqref{eq.P1.5.4b}, the fair outcome ratios, as a function of the probability $p$ for the positive uncertainty outcome $C_{u}$, are adjusted downward, relative to Figure~\ref{fig:risk_seeking_1g}. Furthermore, the fairness symmetry point $p = 0.5$ has been adjusted downward in Figure~\ref{fig:risk_seeking_1i}.

If we have a small initial wealth, Figure~\ref{fig:risk_seeking_1i}, and we stand to gain more than we initially would have gained. Then, for given outcome ratios, we will be more inclined to accept the possibility of gaining nothing, relative to the case where we have a large initial wealth, Figure~\ref{fig:risk_seeking_1g}, as the pay-out, in terms of subjective consequences, is relatively higher under a small initial wealth. 

Comparing Figures~\ref{fig:risk_seeking_1j} and~\ref{fig:risk_seeking_1g}, we see that by taking into account the initial wealth $m$, through the Bernoulli law, \eqref{eq.P1.5.4b}, the fair outcome ratios, as a function of the probability $p$ for the negative uncertainty outcome $C_{u}$, are adjusted upward, relative to Figure~\ref{fig:risk_seeking_1g}. Furthermore, the fairness symmetry point $p = 0.5$ has been adjusted upward in Figure~\ref{fig:risk_seeking_1j}. These adjustments make nothing but sense. 

If we have a small initial wealth, Figure~\ref{fig:risk_seeking_1j}, and we stand to lose more than we initially would have lost. Then, for given outcome ratios, we will be less inclined to accept the possibility of losing even more, relative to the case where we have a large initial wealth, Figure~\ref{fig:risk_seeking_1g}, as the penalty, in terms of subjective consequences, is relatively larger under a small initial wealth.

As our initial wealth tends to infinity, and our utility for money becomes linear, we will perceive both problems to be symmetric, as monetary losses are weighed the same as monetary gains, Figure~\ref{fig:risk_seeking_1g}.   


\section{The Psychological Certainty Effect, Part II}
\label{KT2}
Risk seeking refers to a specific pattern in betting behavior. Uncertain larger gains are preferred over sure smaller gains and uncertain larger losses are preferred over sure smaller losses. The psychologists Kahneman and Tversky state that risk seeking constitutes one of the minimal challenges that must be met by any adequate descriptive theory of choice, \cite{Tversky92}.

The observation that large gains are preferred over sure much smaller gains is commensurate with the fact that we may prefer high-risk, high-yield investment opportunities over low-risk, low-yield ones. Likewise, the observation that uncertain larger losses are preferred over sure smaller, though still substantial, losses is in accordance with those instances in the past where traders incurred hundreds of millions in losses, in their attempts to make good on their previous losses\footnote{As, for example, happened to Nicholas William Leeson, a trader for the Barrings Bank in the nineties. Though we believe that Leeson would have acted less recklessly had he been investing his own money, instead that of the deposit holders. That is, we expect that his Weber constant for his own money, say, $q$, was markedly larger than his Weber constant for the deposit holders money, say, $q_{0}$, where $0 \leq q_{0} <\,< q$.}. 

If the signs of the outcomes in the risk seeking betting scenarios are reversed, then the preferences between the bets will also reverse. This is called the reflection effect, \cite{Kahneman79}. So, risk seeking in the positive domain is accompanied by risk aversion in the negative domain. Conversely, risk seeking in the negative domain is accompanied by risk aversion in the positive domain.

\subsection{Risk Seeking I}  
We first give an example of risk seeking in the case of a small probability of winning a large prize, that is, risk seeking in the positive domain. This case of risk seeking represents our tendency to profit maximization and demonstrates that we will be willing to invest in a long shot if the pay-out is high enough. 

The outcome probability distributions for the respective bets in our risk seeking example are 
\begin{equation}
	\label{eq.P1.10.1}
	 p\!\left(\left.O\right|D_{1}\right) = \begin{cases} 
	 										   0.001, \qquad O = 5000 \\
										     0.999, \qquad O = 0
											\end{cases}
\end{equation} 
and
\begin{equation}
	\label{eq.P1.10.2}
	 p\!\left(\left.O\right|D_{2}\right) = \begin{cases} 
	 											1.0, \qquad O = 5 
											\end{cases}
\end{equation} 
It is found that 72\% of $N = 72$ subjects prefer decision $D_{1}$ over $D_{2}$, \cite{Kahneman79}. Even though both bets have the same expectation value of
\[
	E\left(\left.O\right|D_{1}\right) = 0.001 \times 5000 = 5 = 1.0 \times 5 = E\left(\left.O\right|D_{2} \right) .
\]  
We now interpret this finding in terms of the Bayesian decision theoretic framework.

Seeing that Kahneman and Tversky's Figure~3 of their \cite{Tversky92} is in close correspondence with Figure~\ref{fig:risk_seeking_1g}, we will assume here a linear utility for monetary outcomes\footnote{This assumption is not that far-fetched, seeing that we have here imaginary increments in an imaginary inital wealth; stated differently, no money will be actually gained or lost, as a consequence of our decisions. See also Appendix~\ref{Framing} for a further discussion of the validity of hypothetical bets.}. This assumption allows us to use either Figure~\ref{fig:risk_seeking_1e} or Figure~\ref{fig:risk_seeking_1f} to `read off' the corresponding fair probability. 

It is found that the fair probability for the certainty bets \eqref{eq.P1.10.1} and \eqref{eq.P1.10.2} is given as:
\begin{equation}
	\label{eq.P1.10.3}
	 P_{fair} = 3.985 \times 10^{-6}.
\end{equation}
This fair probability may be checked to correspond with the fair interval, \eqref{eq.P1.10.31}:
\begin{equation}
	\label{eq.P1.10.4}
	\left(0, 10\right) = \left(0, 2 O_{c}\right),
\end{equation}
where the overshoot in the lower bound is re-set to zero\footnote{See also \eqref{eq.why.15b}.}.

So, if the probability of the uncertain gain exceeds this fair probability, as it does, then, as our `expected' utility moves away from the equilibrium situation to a greater expected gain, we will accept the uncertainty bet $D_{1}$. As the gain in the utility upper bound under $D_{1}$ will dominate the gain in the utility lower bound under $D_{2}$. 

And indeed, it is found that 72\% of $N = 72$ subjects prefer decision $D_{1}$ over $D_{2}$, \cite{Kahneman79}, even though both bets have the same outcome expectation values. The phenomenon of utility upper bound dominance for gains constitutes risk seeking in the positive domain.

\subsection{Risk Aversion I}
The above analysis may also be performed for the case when there is a small probability of loosing a large sum of money. We then will see a reversal in the preference for bet $D_{1}$ over bet $D_{2}$ to a preference for bet $D_{2}$ over bet $D_{1}$. Risk aversion in the negative domain represents our tendency to hedge against large and catastrophic losses. 

The outcome probability distributions for the respective bets are\footnote{Compare with \eqref{eq.P1.10.1} and \eqref{eq.P1.10.2}.}:  
\begin{equation}
	\label{eq.P1.10.10b}
	 p\!\left(\left.O\right|D_{1}\right) = \begin{cases} 
	 											0.001, \qquad O = -5000 \\
												0.999, \qquad O = 0
											\end{cases}
\end{equation} 
and
\begin{equation}
	\label{eq.P1.10.10c}
	 p\!\left(\left.O\right|D_{2}\right) = \begin{cases} 
	 											1.0, \qquad O = -5 
											\end{cases}
\end{equation}
It is found that 83\% of $N = 72$ subjects preferred the bet $D_{2}$ over $D_{1}$, \cite{Kahneman79}. 

We will again assume here a linear utility for monetary outcomes. This assumption allows us to use Figure~\ref{fig:risk_seeking_1g} to `read off' the corresponding fair probability. 

It is found that the fair probability for the certainty bets \eqref{eq.P1.10.10b} and \eqref{eq.P1.10.10c} is given as:
\begin{equation}
	\label{eq.P1.10.3b}
	 P_{fair} = 3.985 \times 10^{-6}.
\end{equation}
This fair probability may be checked to correspond with the fair interval, \eqref{eq.P1.10.31b}:
\begin{equation}
	\label{eq.P1.10.4b}
	\left(-10, 0\right) = \left(2 O_{c}, 0\right),
\end{equation}
where the overshoot in the upper bound is re-set to zero\footnote{See also \eqref{eq.why.15b}.}.

So, if the probability of the uncertain loss exceeds this fair probability, as it does, then as our `expected' utility moves away from the equilibrium situation to a greater expected loss, we will reject the uncertainty bet $D_{1}$. As the gain in the utility lower bound under $D_{2}$ will dominate the gain in the utility upper bound under $D_{1}$. 

And indeed, it is found that 83\% of $N = 72$ subjects prefer decision $D_{2}$ over $D_{1}$, \cite{Kahneman79}, even though both bets have the same outcome expectation values. The phenomenon of utility lower bound dominance for losses constitutes risk aversion in the negative domain.

\subsection{Risk Seeking II}  
We now give an example of risk seeking when people must choose between a sure loss and a substantial probability of a larger loss, that is, risk seeking in the negative domain. This case of risk seeking represents our tendency to try to evade large and catastrophic losses. 

The outcome probability distributions for the respective bets in our risk seeking example are 
\begin{equation}
	\label{eq.P1.10.11}
	 p\!\left(\left.O\right|D_{1}\right) = \begin{cases} 
												0.8, \qquad O = -4000 \\
												0.2, \qquad O = 0
											\end{cases}
\end{equation} 
and
\begin{equation}
	\label{eq.P1.10.12}
	 p\!\left(\left.O\right|D_{2}\right) = \begin{cases} 
												1.0, \qquad O = -3000 
											\end{cases}
\end{equation} 
It is found that 92\% of $N = 95$ subjects preferred the bet $D_{1}$ over $D_{2}$, \cite{Kahneman79}. 

If we assume a linear utility for monetary outcomes, then we may use Figure~\ref{fig:risk_seeking_1g} to `read off' the corresponding fair probability. It is found that the fair probability for the certainty bets \eqref{eq.P1.10.11} and \eqref{eq.P1.10.12} is given as:
\begin{equation}
	\label{eq.P1.10.3c}
	 P_{fair} = 0.8536.
\end{equation}
This fair probability may be checked to correspond with the fair interval, \eqref{eq.P1.10.32b}:
\begin{equation}
	\label{eq.P1.10.4c}
	\left(-4000, -2000\right) = \left(O_{u}, 2 O_{c} - O_{u}\right),
\end{equation}
where the overshoot in the lower bound is re-set to minus four thousand\footnote{See also \eqref{eq.why.15b}.}.

So, if the probability of the uncertain loss is less than the fair probability \eqref{eq.P1.10.3c}, as it is, then as our `expected' utility moves away from the equilibrium situation to a lesser expected loss, we will accept the uncertainty bet $D_{1}$. As the gain in the utility upper bound under $D_{1}$ will dominate the gain in the lower bound under $D_{2}$. 

And indeed, it is found that 92\% of $N = 95$ subjects prefer decision $D_{1}$ over $D_{2}$, \cite{Kahneman79}, even though $D_{1}$ has a slightly larger outcome expectation value than $D_{2}$. The phenomenon of utility upper bound dominance for losses constitutes risk seeking in the negative domain.

\subsection{Risk Aversion II}
The previous analysis may also be performed for the opposite case of a sure gain and a substantial probability of a larger gain. We then will see a reversal in the preference for bet $D_{1}$ over bet $D_{2}$ to a preference for bet $D_{2}$ over bet $D_{1}$. Risk aversion in the positive domain represents our tendency to secure our profits. 

The outcome probability distributions for this problem of choice are\footnote{Compare with \eqref{eq.P1.10.11} and \eqref{eq.P1.10.12}.}:
\begin{equation}
	\label{eq.P1.10.24}
	 p\!\left(\left.O\right|D_{1}\right) = \begin{cases} 
												0.8, \qquad O = 4000 \\
												0.2, \qquad O = 0
											\end{cases}
\end{equation} 
and
\begin{equation}
	\label{eq.P1.10.25}
	 p\!\left(\left.O\right|D_{2}\right) = \begin{cases} 
												1.0, \qquad O = 3000 
											\end{cases}
\end{equation} 
It is found that 80\% of $N = 95$ subjects preferred the bet $D_{2}$ over $D_{1}$, \cite{Kahneman79}. 

If we assume a linear utility for monetary outcomes, then we may use either Figure~\ref{fig:risk_seeking_1e} or Figure~\ref{fig:risk_seeking_1f} to `read off' the corresponding fair probability. It is found that the fair probability for the certainty bets \eqref{eq.P1.10.24} and \eqref{eq.P1.10.25} is given as:
\begin{equation}
	\label{eq.P1.10.3d}
	 P_{fair} = 0.8536.
\end{equation}
This fair probability may be checked to correspond with the fair interval, \eqref{eq.P1.10.32}:
\begin{equation}
	\label{eq.P1.10.4d}
	\left(2000, 4000\right) = \left(2 O_{c} - O_{u}, O_{u}\right),
\end{equation}
where the overshoot in the upper bound is re-set to four thousand\footnote{See also \eqref{eq.why.15b}.}.

So, if the probability of the uncertain loss is less than the fair probability \eqref{eq.P1.10.3d}, as it is, then as our `expected' utility moves away from the equilibrium situation to a lesser expected gain, we will accept the certainty bet $D_{2}$. As the gain in the utility lower bound under $D_{2}$ will dominate the gain in the upper bound under $D_{1}$. 

And indeed, it is found that 80\% of $N = 95$ subjects prefer decision $D_{2}$ over $D_{1}$, \cite{Kahneman79}, even though $D_{2}$ has a slightly larger outcome expectation value than $D_{1}$. The phenomenon of utility lower bound dominance for gains constitutes risk aversion in the positive domain.

\section{Discussion}
In this fact sheet we have presented the case for the Bayesian decision theory. It may be read in Jaynes' \cite{Jaynes03}, that to the best of his knowledge, there are as of yet no formal principles at all for assigning numerical values to loss functions; not even when the criterion is purely economic, because the utility of money remains ill-defined. In the absence of these formal principles, Jaynes final verdict was that decision theory can not be fundamental. This situation may have changed. 

The Bernoulli utility function, initially derived by Bernoulli, by way of common sense first principles\footnote{See Appendix~\ref{BernoulliA}.}, has now been derived by way of a consistency argument\footnote{See Section~\ref{BernoulliC}.}. This might explain why it is that Bernoulli's utility function has proven to be so ubiquitous and successful the field of sensory perception research; simply because it is consistent. 

The first two algorithmic steps of the Bayesian decision theory, respectively, the construction of outcome probability distributions by way of the Bayesian probability theory and the construction of utility probability distributions by way of the Bernoulli utility function, allow us no freedom. To construct our outcome and utility probability distributions otherwise, would be to invite inconsistency. 

Our initial justification for the Bernoulli utility function had come from the observation that this function, in the guise of the Weber-Fechner and the Steven's power law, had been demonstrated by psycho-physics to be an appropriate model for the way we humans perceive the increments in sensory stimuli, in terms of sensation strength. So, if monetary outcomes are considered to be a sensory stimuli, in the most abstract sense of the word, then it would follow the Bernoulli utility function would be the most appropriate model for the way we humans perceive the increments in monetary wealth.

As in the course of our research we came to trust the Bayesian decision algorithm to teach our intuition, in those instances where the intuitive `resolution' is lacking to make clear and crisp choices\footnote{Just like we have learned, having been Bayesians for the past ten years, to trust the Bayesian probability algorithm to teach our intuition, in those instances where the intuitive resolution is lacking to make clear and crisp plausibility assessments; that is, by analogy, Jaynes' reasoning computer of Bayesian probability theory, \cite{Jaynes03}, had become a decision making computer.}, the burden fell on us to provide a proof of the fundamentalness of the Bayesian decision theory. 

The history of Bayesian probability theory has taught us that the usefulness of a theory, in terms of its practical and beautifully intuitive results, in the absence of a compelling axiomatic basis, provides no safeguard against attacks by those who choose to close their eyes to this usefulness\footnote{Note that this historical fact explains why Bayesians have their axiomatic house in such good order. This process started with the work of Cox, \cite{Cox46}, was expanded upon by Jaynes, \cite{Jaynes03}, which was then further refined by the work of Knuth and Skilling, \cite{Knuth10}. Moreover, the more general axiomatic framework of the latter has enabled them, amongst other things, \cite{Knuth14}, to bring some order to the field of quantum theory, by showing why this theory is forced to use a complex arithmetic, \cite{Goyal10}.}. This is why we felt compelled to search for a consistency derivation of the Bernoulli utility function. 

Especially so, since we had taken painstaking care to search out those `unyielding practical realities', that would put our foundations to the test. And it had been found that all these practical realities fell nicely in line with the proposed foundations of the Bayesian decision theory, as may be witnessed in our treatment of the Ellsberg paradox and the Kahneman and Tversky data on the certainty effect\footnote{The Allais paradox corresponds, more or less, with Risk Aversion~II in Section\ref{KT2}. But one may take, trivially, any of the proposed Allais paradoxes and solve them with the here proposed decision algorithm.}

So, having presented a consistency proof for the Bernoulli utility function, the question now is: Is the Bayesian decision theory, just like the Bayesian probability and information theories, Bayesian in the strictest sense in the word, or, equivalently, an inescapable consequence of the desideratum of consistency? We will now try to answer this question.

There is one degree of freedom remaining in the Bayesian decision theory as a whole. This remaining degree of freedom is the criterion of choice which states that we should maximize the sum of the lower bound of the utility probability distribution\footnote{See Section~\ref{Choice}.}. 

In any problem of choice we will endeavor to choose that decision which has a corresponding utility probability distribution that is lying most the right on the utility axis; that is, we will choose to maximize our utility probability distributions. In this there is little freedom. But we are free, in principle, to maximize the positions of our utility probability distributions any way we see fit. Nonetheless, we believe that it is always a good policy to take into account all the pertinent information we have. 

For example, if we only maximizes the expected value of the utility probability distribution, then we will, by definition, neglect the information that the standard deviation of the utility probability distribution has to bear on our problem of choice, by way of the symmetry breaking in the case of an overshoot of one of the bounds. Likewise, we are free to only maximize one of the bounds of our utility probability distributions, while neglecting the other. But in doing so, we will neglect the possibility of either (catastrophic) losses in the lower bound or (astronomical) gains in the upper bound. 

In the Bayesian probability theory we have an analogous situation of both constraint and freedom. The Bayesian probability theory states that the product and sum rules of probability theory are the only two consistent and, therefore, admissible operators to combine probabilities with. But it also states that in our assignment of our probabilities we are totally free to do as we see fit. There are no `wrong' probability distributions, only better or worse informed ones. It's all a matter of choice.
\\
\\
\noindent\textbf{Acknowledgments:} The research leading to these results has received partial funding from the European Commission's Seventh Framework Program [FP7/2007-2013] under grant agreement no.265138.

\appendix

\section{Bayesian Probability Theory}
\label{Bayes}
The whole of Bayesian probability theory flows forth from two simple rules. The product rule,
\begin{equation}
	\label{eq.product}
	P\!\left(A\right) P\!\left(\left.B\right|A\right) =	P\!\left(A B\right) = P\!\left(B\right) P\!\left(\left.A\right|B\right), 
\end{equation}
and the sum rule
\begin{equation}
	\label{eq.sum}
	P\!\left(\overline{A}\right) = 1- P\!\left(A\right).
\end{equation}

By way of the product and the sum rule, we may derive the generalized sum rule,
\begin{equation}
	\label{eq.genSum}
	P\!\left(A + B\right) = P\!\left(A\right) + P\!\left(B\right) - P\!\left(A B\right).
\end{equation}
If we have that the propositions are exhaustive and mutually exclusive, that is, $B = \overline{A}$, we may derive, by way of \eqref{eq.sum}, \eqref{eq.genSum}, and the fact that $p\!\left(A \overline{A}\right) = 0$, the most primitive probability distribution:
\begin{equation}
	\label{eq.genSum2a}
	P\!\left(A + \overline{A}\right) = P\!\left(A\right) + P\!\left(\overline{A}\right) = 1,
\end{equation}
This probability distribution then may be further generalized, by taken as it propositional elements the exhaustive and mutual exclusive conjunctions $A_{i} B_{j}$, to the bivariate probability distribution:
\begin{equation}
	\label{eq.genSum2}
	\sum_{i} \sum_{j} P\!\left(A_{i} B_{j}\right) = 1,
\end{equation}
which allows us to `marginalize' over the parameter, say, $B_{j}$, which is of no direct interest:
\begin{equation}
	\label{eq.genSum2c}
	P\!\left(A_{i} \right) = \sum_{j} P\!\left(A_{i} B_{j}\right),
\end{equation}
where
\begin{equation}
	\label{eq.genSum2d}
	\sum_{i}  P\!\left(A_{i} \right) = 1.
\end{equation}

We may let $i$ and $j$, that is, the number of propositions $A_{i}$ and $B_{j}$ in \eqref{eq.genSum2}, tend to infinity. By doing so, we go from discrete to continuous probability distributions. Furthermore, we may add propositions $C_{k}$, $D_{l}$, $E_{m}$, etc..., and so get higher variate distributions.

Now, to a non-Bayesian it may seem to be somewhat surprising, that the whole of Bayesian probability theory flows forth from the product and rules. But the whole of Boolean logic\footnote{We use the term `Boolean algebra' in its meaning as referring to two-valued logic in which symbols like `A' stand for propositions, \cite{Jaynes03}.}, on an operational level, is also captured by the AND- and NOT-operations. These operations correspond, respectively, with \eqref{eq.product} and \eqref{eq.sum}; as these operators combine, with the negation of a NAND-operation, in the OR-operation, which corresponds with \eqref{eq.genSum}. 

Moreover, it may be shown that Boolean logic is just a special limit case of the more general Bayesian probability theory. The operators of Boolean logic combine in a like manner as the operators in Bayesian probability theory. But in Boolean logic propositions can have only the truth values true or false. Whereas in Bayesian probability theory propositions can have plausibility values in the interval $\left[0, 1\right]$, where $0$ and $1$, respectively, correspond with false and true. 

So Boolean logic is the language of deduction, whereas Bayesian probability theory is the language of both induction and deduction; the former being a limit case of the latter, in which we have absolute knowledge about the propositions in play. 

Now, on the conceptual level Bayesian probability theory is very simple. However, on an implementation level, when doing an actual data-analysis, it may be quite challenging\footnote{In close analogy, Boolean logic, which is a specific limit case of Bayesian probability theory, is simple on the conceptional level. However, on the implementation level it may be quite challenging, when, say, we use this Boolean logic to design logic circuits for computers.}; which, as an aside, makes it fun to do Bayesian statistics. And we refer the interested reader to \cite{Skilling05}, for a first cursory overview on the considerations that come with a Bayesian data-analysis\footnote{Though the absolute autority is \cite{Jaynes03}. But the reading of this 680-page tome would require a considerable time investment on the part of the reader. But then again, as Calculus is the royal highway to the exact sciences, so we have that Jaynes' \textit{Probability Theory: The Logic of Science} is the royal highway to Bayesian statistics.}.
 
\section{The Ubiquitous Bernoulli Utility Function}
\label{BernoulliA}
The utility of a given outcome is the perceived worth of that outcome. If we take the utilities that monetary outcomes hold for us to be an incentive for our decisions, then we may perceive money to be a stimulus. 

For the rich man ten dollars is an insignificant amount of money. So, the prospect of gaining or losing ten dollars will fail to move the rich man, that is, an increment of ten dollars for him has an utility which tends to zero. 

For the poor man ten dollars is two days worth of groceries and, thus, a significant amount of money. So, the prospect of gaining or losing ten dollars will most likely move the poor man to action. It follows that an increment of ten dollars for him has an utility significantly greater than zero.

We now will give the derivations of the Bernoulli, the Weber-Fechner, and Steven's power laws. It will be seen all that these three laws are equivalent.

\subsection{The Bernoulli utility function.}
Consider persons $A$ and $B$, with $A$ having a fortune of 100.000 full-ducats, and with $B$ a fortune of 100.000 semi-ducats, a semi-ducat being the half of a full-ducat. Let $f_{A}$ and $f_{B}$ be the moral value functions, defined on, respectively, the monetary full-ducat axis $x$ and the semi-ducat axis $\tilde{x}$. Let $x_{A}$ and $\tilde{x}_{B}$ stand for the initial wealths of $A$ and $B$, respectively;  where $x_{A}$ and $\tilde{x}_{B}$ are points on the monetary axes $x$ and $\tilde{x}$, respectively.

Bernoulli derived his law by way of three simple variance considerations for the moral functions $f_{A}$ and $f_{B}$, \cite{Bernoulli38, Masin09}:
\begin{enumerate}
	\item For an arbitrary increment $c$ in wealth, the moral movement of this increment will be less for the rich man, than for the poor man; that is, if we make for $f_{B}$ the appropriate change of variable, from $\tilde{x}$ to $x$, then we have that 
	\[
	\left.\frac{d}{d x} f_{A}\left(x\right)\right|_{c} < \left.\frac{d}{d x} f_{B}\left(x\right)\right|_{c}.
	\]
	From which it follows that effect of $c$ on a given $f$ decreases as the initial wealth increases. \\
	\item It is proposed that the movement in a general moral value function $f$, for a given positive increment $dx$, is proportional to the value of this increment; that is, 
	\[
		\left.\frac{d}{d u} f\left(u\right)\right|_{c = d x} \propto d x,
	\]
	as this is the simplest function for which $f$ increases as a function of an increment in $x$. \\
	\item Furthermore, it is proposed that this movement in $f$ is inversely proportional to the value of the initial wealth $x$; that is, 
	\[
			\left.\frac{d}{d u} f\left(u\right)\right|_{c = dx} \propto \frac{1}{x}.
	\]
\end{enumerate}
where `$\propto$' is the proportionality sign.

Bernoulli arrived at his third consideration, using the following reasoning. The change in moral value of $c$ full-ducats for $A$ will be half the change in moral value of $c$ full-ducats for $B$. Only if either $B$ sees his fortune increased to 200.000 semi-ducats, or, equivalently, 100.000 full-ducats, or if $A$ sees his fortune reduced to 50.000 full-ducats, or, equivalently, 100.000 semi-ducats, only then will $B$ have the same change in moral value as $A$ for $c$ full-ducats\footnote{In these two limit cases the Bernoulli variance argument will tend to the invariance argument \eqref{BerC.1}, which was used in Section~\ref{BernoulliC}.}. So, if we make for $f_{B}$ the appropriate change of variable from $\tilde{x}$ to $x$, 
\begin{equation}
 \label{eq.Ber.1}
  \frac{\left.\frac{d}{d x} f_{A}\left(x\right)\right|_{c}}{\left.\frac{d}{d x} f_{B}\left(x\right)\right|_{c}} =  \frac{1}{2} = \frac{x_{B}}{x_{A}}, 
\end{equation}
where $x_{B}$ is the initial fortune of $B$, translated from the semi-ducat $\tilde{x}$-axis to the full-ducat $x$-axis.  

It follows from \eqref{eq.Ber.1} that we have, in general, that the change in moral value is inversely proportional to the initial we hold, that is,
\begin{equation}
 \label{eq.Ber.2}
   \left.\frac{d}{d x} f\left(x\right)\right|_{c} \propto  \frac{1}{x}, 
\end{equation}
which is Bernoulli's third consideration.

If we combine the second and the third consideration, we obtain the differential equation
\begin{equation}
 \label{eq.Ber.3}
   f'\left(x\right) =  q \frac{dx}{x}, 
\end{equation}
which, if solved for the boundary condition that for a given person with an initial wealth of $x_{0}$ an increment of zero holds no utility, either negative or positive, gives
\begin{equation}
 \label{eq.Ber.4a}
   f\left(x\right) =  q \log \frac{x}{x_{0}}, 
\end{equation}  
which may be rewritten as
\begin{equation}
 \label{eq.Ber.4}
   f\left(\left.\Delta x\right|x_{0}\right) =  q \log \frac{x_{0} + \Delta x}{x_{0}}. 
\end{equation} 

\subsection{The Weber-Fechner Law.}
Let $S$ signify stimuli intensity and let $Q$ signify sensation strength. Weber's law states that the increment $\Delta S$ needed to elicit a judgment that $S + \Delta S$ is just noticeably different from $S$ is proportional to $S$:  
\begin{equation}
	\label{eq.P1.8.1}
	  \Delta S = w S,
\end{equation}    
where $w$ is a positive constant dependent upon the specific type of sensory stimulus offered and $\Delta S$ is understood to be the stimulus increment corresponding with a just noticeable difference. 

Fechner generalized the experimental Weber law by stating that all differences in sensational strength, and not only the ones that are just noticeable, are proportional to the relative change $\Delta S/S$, that is,
\begin{equation}
	\label{eq.P1.8.2}
	 \Delta Q = q \frac{\Delta S}{S}.
\end{equation}    
where $k$ is a positive constant dependent upon the specific type of sensory stimulus offered and $\Delta S$ is now understood to be the stimulus increment corresponding with the increment in sensation strength $\Delta Q$. 

Dividing both sides of \eqref{eq.P1.8.2} by $\Delta S$ gives
\begin{equation}
	\label{eq.P1.8.3}
	 \frac{\Delta Q}{\Delta S} = q \frac{1}{S}.
\end{equation}    
Fechner then makes the assumption that, just as a physically small quantity $\Delta S$ can be reduced without limit to the differential $dS$, so a small quantity of sensation can be reduced without limit to the differential $dQ$. By way of this assumption, we may let \eqref{eq.P1.8.3} tend to the differential equation
\begin{equation}
	\label{eq.P1.8.4}
	 \frac{dQ}{dS} = q \frac{1}{S}.
\end{equation}   
The general solution of this differential equation is 
\begin{equation}
	\label{eq.P1.8.5}
	Q = q \log S + c,
\end{equation} 
where $c$ is some constant of integration. 

Introducing an initial value condition for \eqref{eq.P1.8.4} that says that at stimulus value $S_{0}$ there is no sensation strength, that is, $Q\left(S_{0}\right) = 0$, leaves us with the Weber-Fechner law  
\begin{equation}
	\label{eq.P1.8.6}
	Q\left(\left.S\right|S_{0}\right) = q \log \frac{S}{S_{0}},
\end{equation}
or, equivalently, 
\begin{equation}
	\label{eq.P1.8.7}
	Q\left(\left.\Delta S\right|S_{0}\right) = q \log \frac{S_{0} + \Delta S}{S_{0}}.
\end{equation}

The Weber-Fechner law, \eqref{eq.P1.8.7}, is identical to the utility function which had been proposed a century earlier by Bernoulli,  \eqref{eq.Ber.4}. 

Fechner himself was aware of this equivalence. Nonetheless, he believed his derivation to be the more general. Fechner argued that Bernoulli's derivation only applied to the special case of utility, whereas his law, though identical, applied to all sensations, as it invokes Weber's law. 

However, as pointed out in \cite{Masin09}, Fechner failed to provide any compelling reason why the principles employed in Bernoulli's derivation of the subjective value of objective monies should not be extendible to sensations in general. Nonetheless, we do believe that Fechner acted in good faith, in denying Bernoulli scientific primacy. 

First of all, Fechner called the Weber-Fechner law, when he first published it, very modestly, the Weber law. Second of all, Fechner had a deep spiritual need for some kind of harmony between the physical and mental universes, and the Weber-Fechner law provided him with this harmony, for this law spoke of the basic oneness of the physical and mental universes, \cite{Fancher90}.

The Weber-Fechner law demonstrated that both universes adhered to seemingly mechanistic laws. It then followed that the freedom of the latter universe, in terms of free will and volition, implied, by way of analogy, a commensurate freedom of the former; thus, opening the way for the possibility of a besouled physical universe. Which had become Fechner's only hope for spiritual salvation, \cite{Fancher90}. 

We can imagine that Fechner might have felt that a law that assigned subjective values to objective monies was too arbitrary and sordid a foundation for the lofty purpose he wished it to serve. In contrast, the initial Weber law allowed Fechner to forgo of the money argument and derive a law, which though in form identical to Bernoulli's, differed in that it applied to all human sensations. 

\subsection{Steven's Power Law.}
Steven's power law is based on the observation, that it is the ratio $\Delta Q/Q$, rather than the difference $\Delta Q$, that is proportional to  $\Delta S/S$, \cite{Stevens61}. This observation leads to the equality
\begin{equation}
	\label{eq.P1.8.8}
	 \frac{\Delta Q}{Q} = q \frac{\Delta S}{S}.
\end{equation}    
Letting the differences in $Q$ and $S$ go to differentials, we may rewrite \eqref{eq.P1.8.8} as
\begin{equation}
	\label{eq.P1.8.9}
	 \frac{d Q}{Q} = q \frac{d S}{S}.
\end{equation} 
This equation has its general solution 
\begin{equation}
	\label{eq.P1.8.10}
	 \log Q = q \log S + c'.
\end{equation}
Taking the exponent of both sides of \eqref{eq.P1.8.10}, we get the power law for stimulus perception 
\begin{equation}
	\label{eq.P1.8.11}
	 Q = c S^{q},
\end{equation}
where $c = \exp\left(c'\right)$.

Stevens found the power law to hold for several sensations; binaural and monaural loudness, brightness, lightness, smell, taste, temperature, vibration duration, repetition rate, finger span, pressure on palm, heaviness, force of hand grip, autophonic response, and electric shock, \cite{Stevens61}. 

The power law is applied by letting subjects compare the sensation ratio of $Q_{1}$ to $Q_{0}$ for corresponding stimuli strengths $S_{1}$ and $S_{0}$:
\begin{equation}
	\label{eq.P1.8.12}
	 \frac{Q_{1}}{Q_{0}} = \left(\frac{S_{1}}{S_{0}}\right)^{q}.
\end{equation}
Let $S_{1} = S_{0} + \Delta S$, where $\Delta S$ is some increment, then we may rewrite \eqref{eq.P1.8.12} as
\begin{equation}
	\label{eq.P1.8.13}
	 \frac{Q_{1}}{Q_{0}} = \left(\frac{S_{0} + \Delta S}{S_{0}}\right)^{q}.
\end{equation}
For an increment of $\Delta S = 0$, the ratio of perception stimuli will be $Q_{1}/Q_{0} = 1$. Taking the log of the ratio \eqref{eq.P1.8.13} we may map the ratio of perceived stimuli to a corresponding utility scale where a zero increment $\Delta S$ corresponds with a zero utility:
\begin{equation}
	\label{eq.P1.8.14}
	 Q'\left(\left.\Delta S\right|S_{0}\right) = \log \frac{Q_{1}}{Q_{0}} = q \log \frac{S_{0} + \Delta S}{S_{0}}
\end{equation} 
But this is just the Weber-Fechner law, \eqref{eq.P1.8.8}. 

\subsection{Summary.} The Weber-Fechner law gives us just noticeable differences on a log scale, \eqref{eq.P1.8.8}. The power law gives us ratios of sensation strengths, \eqref{eq.P1.8.13}. Taking the log of the ratio of sensation strengths, we may obtain the just noticeable differences again, \eqref{eq.P1.8.14}. But the Weber-Fechner for just noticeable differences is just the Bernoulli utility function for utilities, \eqref{eq.Ber.4}. 

We refer the reader to \cite{Masin09}, for a discussion of Thurnstone's derivation of the satisfaction law. This law, which takes as its input the increment in the number of items of commodity, is also of the form of Bernoulli's utility function.

\section{The Negative Bernoulli utility function}
\label{BernoulliB}
In this appendix we present the negative Bernoulli utility function for debts, which is a corollary of the Bernoulli utility function for income. The negative Bernoulli utility function predicts that for the very poor, having a small initial wealth and large initial debts, a large loss of direct income will be more devastating, than an increase of, say, twice that loss in their long-term debt. This law also explains why, for these poor, having a small initial wealth and large initial debts, the temptation to take out loans, if offered the opportunity, will be quite great, \cite{Hudson06}.    

Until now we have treated only the case were the maximal loss did not exceed the initial wealth $m$. However, in real life we may lose more than we actually have, by way of debt. So, we now proceed to assign utilities to increments in debt. 

According to the Weber-Fechner law we cannot lose more money than we initially had. Otherwise we may have that the ratio in the logarithm in the Weber-Fechner Law, \eqref{eq.P1.5.3b},
\begin{equation}
	\label{eq.P1.zq.1}
	u\!\left(\left.\Delta S\right|S\right) = q \log \frac{S + \Delta S}{S},
\end{equation} 
may become negative, leading to a breakdown of the logarithm. 

However, whenever we incur a debt we lose more money than we have. Furthermore, we can have a debt and an income, both at same time. So, we propose that there are two different monetary stimuli dimensions in play; the first dimension being an actual income dimension and the second dimension being a debt dimension. 

We propose to model the debt utilities by way of the negative Weber-Fechner law:
\begin{equation}
	\label{eq.P1.zq.2}
	u\!\left(\left.\Delta D\right|D\right) = - b \log \frac{D + \Delta D}{D}, 
\end{equation} 
where we let $D$ be the initial debt, $\Delta D$ the increment in debt, and $b$ the the Weber constant of a monetary debt.  

The rationale behind \eqref{eq.P1.zq.2} is as follows. If we view a debt increment as a stimulus, then it follows that we may use the psycho-physical Weber-Fechner law in the determination of the moral value of a given debt increment. 

For positive increments $\Delta D$, there is an increase in current debt, whereas for negative increments $\Delta D$, there is a decrease in current debt. In order to assign both a negative utility to an increase in current debt and positive utility to a decrease in debt, we need to multiply the Weber-Fechner law times minus one, \eqref{eq.P1.zq.2}.

If we have no initial debt, that is, $D = 0$, then \eqref{eq.P1.zq.2} tells us that any positive increment in debt $\Delta D$ would have a utility of minus infinity. This is clearly not realistic. So, in order to model an increment in debt for those who are without debt, we must introduce a minimum significant amount of debt which is equal to  minimum significant amount of income, $\gamma$. 

The threshold amount of debt, $\gamma$, may also be used in the case of $\Delta D = - D$, in order to prevent an infinite utility being assigned to a full repaying of one's debts. Using the concept of the minimum significant amount of debt stimulus, we may modify \eqref{eq.P1.zq.2} as 
\begin{equation}
	\label{eq.P1.zq.10}
	u\!\left(\left.\Delta D\right|D\right) = - b \log \frac{D + \Delta D}{D}, \qquad			 -D + \gamma < \Delta D < \infty,
\end{equation} 

If we want to give a graphical representation of \eqref{eq.P1.zq.10}, then the Weber constant $b$, must be set to some numerical value. 

Say, we have a total debt of forty thousand dollars, in the form of a student loan, which we eventually will have to pay back, but not right now. Then introspection would suggest that a increment or decrement of an amount less than a thousand dollars would not move us that much. 

So, $\Delta D = 1000$ constitutes one utile, or, equivalently, a just noticeable difference in debt for an initial debt of $D = 40.000$, that is, \eqref{eq.P1.5.4b}:
\begin{equation}
	\label{eq.P1.zq.3}
	1 \:\text{utile} = - b \log \frac{40.000 - 1000}{40.000}.
\end{equation} 
If we then solve for the unknown Weber constant $b$ of debt stimuli,
\begin{equation}
	\label{eq.P1.zq.4}
		b = - \frac{1}{\log 390000 - \log 40000}\approx 40,
\end{equation} 
we find this Weber constant to be smaller by a factor of $2.5$ than the Weber constant $q$ of income stimuli, \eqref{eq.P1.10.2b}. 

It is well possible that this difference in Weber constants can be attributed to the difference in abstractness of the concepts. The losing of actual monies is quite concrete, whereas the accrueing of a debt, repayable somewhere in a distant future, is somewhat more abstract. 

But there is always a chance that these authors were off in their introspection\footnote{Note that actual value of the Weber constant $b$ of the debt stimuli has no direct bearing on any of the results given in this fact sheet; save the handful of examples which are given in this section, in order to demonstrate the qualitative behavior of the negative Bernoulli utility function, or, equivalently, the negative Weber-Fechner law.}, and that both Weber constants should be approximately equal. We leave this issue, together with the psychological reality of the phenomenon of debt relief, given below, for future psychological experimentation, as we proceed with our discussion of the debt utilities. 

\begin{figure}[h]
		\centering
			\includegraphics[width=0.50\textwidth]{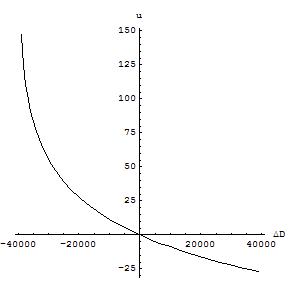}
		\caption{Utility plot for initial debt 40.000 dollars}
		\label{fig:P1.zq.1}
	\end{figure}

Suppose that a student has a student loan which has accumulated to forty thousand dollars. Then, by way of \eqref{eq.P1.zq.2} and \eqref{eq.P1.zq.4}, we obtain the following mapping of increments in debt to utilities, Figure~\ref{fig:P1.zq.1}.

As stated previously, loss aversion is the phenomenon that losses may loom larger than gains. In Figure~\ref{fig:P1.zq.1} we see the phenomenon that debt reduction may loom larger than debt increase. We will call this corollary of the psycho-physical Weber-Fechner law: `debt relief', the relief of loosing one's debts. 

Now, does the phenomenon of debt relief correspond with a real psychological phenomenon? We belief that it actually does. 

Say, we have a debt of a a thousand dollars. Then we can imagine ourselves feeling greatly reliefed, were we to be released of our debt. Now, were our debt, instead, to be doubled to two thousand dollars, then we can also imagine ourselves feeling unhappy about this. But this feeling of unhappiness about the doubling of our debt would be of a lesser intensity than the corresponding relief of having our debt acquitted. 

We will now look at the practical implications of the negative Bernoulli utility function, \eqref{eq.P1.zq.2}, and its Weber constant $b$, \eqref{eq.P1.zq.4}. 

A student loan initially represents a gain in debt stimulus. This debt makes itself felt, in terms of actual loss of income, only after graduation, the moment the monthly payments have to be paid and take a considerable chunk out of one's actual income. 

Say, that the student of Figure~\ref{fig:P1.zq.1}, having become a PhD, and having a net income of fifteen hundred dollars, is called upon to make good on his loan, by way of monthly payments of five hundred dollars. Then these payments represent both a loss in income, having a negative utility of, \eqref{eq.P1.zq.1} and \eqref{eq.P1.10.2b}:
\begin{equation}
	\label{eq.P1.zq.5}
	u^{\left(\text{income-loss}\right)}  = 100 \log \frac{1500 - 500}{1500} = -41.5, 
\end{equation} 
as well as a decrements in debt, having a positive utility of, \eqref{eq.P1.zq.2} and \eqref{eq.P1.zq.4}:
\begin{equation}
	\label{eq.P1.zq.6}
	u^{\left(\text{debt-decrease}\right)}  = - 40 \log \frac{40000 - 500}{40000} = 0.5, 
\end{equation} 
It follows that our PhD can find little to no comfort in the fact that he is paying of his debt, as he acutely feels the sting of loss of income. This is, together with the difference in Weber constants\footnote{A difference which accounts only for a factor of $2.5$ in the observed differences of the utilities \eqref{eq.P1.zq.5} and \eqref{eq.P1.zq.6}.}, \eqref{eq.P1.10.2b} and \eqref{eq.P1.zq.4}, reflective of the fact that his utility function for income is highly non-linear in the neighborhood of the increment, whereas his utility function for debt is highly linear in that region.

Now, say that we have another PhD, who during his student days lived a more frugal life style and, consequently, only has a debt of two thousand dollars. For this PhD student, when called upon to make good on the loan, the loss of income will be felt just as keenly, with a negative utility of $u = -40.5$, \eqref{eq.P1.zq.5}. However, he will find more satisfaction in the fact that he is paying of his debts, \eqref{eq.P1.zq.2} and \eqref{eq.P1.zq.4}:  
\begin{equation}
	\label{eq.P1.zq.7}
	u^{\left(\text{debt-decrease}\right)}  = - 40 \log \frac{2000 - 500}{2000} = 11.5, 
\end{equation}  
seeing that he has a more curved utility function for debt than our previous PhD student. 

Nonetheless, the first PhD student may feel, after a couple of years of monthly repayments, when his loan has been reduced to twenty thousand dollars, for the first time, as if he has an actual stake in the repayment of his debt, \eqref{eq.P1.zq.2} and \eqref{eq.P1.zq.4}:    
\begin{equation}
	\label{eq.P1.zq.8}
	u^{\left(\text{debt-decrease}\right)}  = - 40 \log \frac{20000 - 500}{20000} = 1.0, 
\end{equation}  
as his debt repayment utility crosses the threshold of the just noticeable difference.

The negative Bernoulli utility function also gives an explanation why for the very poor, having a minimum monthly wage of seven hundred euros, and already having a large debt of, say, twenty thousand euros, a loss of income of, say, five hundred euros, is perceived to be so much more devastating than an increase in debt of, say, a thousand euros. As for this poor person, the loss of actual income has a negative utility of $-125$ utiles and the gain of an increase of has a negative utility of only $-2$ utiles\footnote{Even if $b = 100$, \eqref{eq.P1.zq.4}, this negative utility will only be $-5$ utiles.}.

Likewise, the temptation for the very poor, if offered the opportunity, to take out a loan of a thousand euros will be quite great. As for this poor person, the immediate gain of a direct increase of a thousand euros in income will have a positive utility of $+89$ utiles, whereas the negative utility of an increase in debt of a thousand dollars will have a negative utility of only $-2$ utiles\footnote{Idem.}, \cite{Hudson06}.

\section{Bayesian Inference}
\label{Fred}
Bayesian statistics is not only said to be common sense quantified, but also common sense amplified\footnote{If Bayesian inference were not common sense amplified, then it could not ever hope to enjoy the successes it currently enjoys in the various fields of science; astronomy, astrophysics, chemistry, image recognition, etc...}, having a much higher `probability resolution' than our human brains can ever hope to achieve \cite{Jaynes03}. 

This statement is in accordance with the Kahneman and Tversky finding that, if presented with some chance of a success, $p$, subjects fail to draw the appropriate binomial probability distribution of the number of successes, $r$, in $n$ draws. Even though subjects manage to find the expected number of successes, they fail to accurately determine the probability spread of the $r$ successes. Kahneman and Tversky see this as evidence that humans are fundamentally non-Bayesian in the way they do their inference, \cite{Kahneman72}. 

We instead propose that human common sense is not hard-wired for problems involving sampling distributions. Otherwise there would be no need for such a thing as data-analysis, as we only would have to take a quick look at our sufficient statistics after which we then would draw the probability distributions of interest. 

However, humans do seem to be hard-wired for more `Darwinian' problems of inference. 

For example, if we are told that our burglary alarm has gone off, after which we are also told that a small earthquake has occurred in the vicinity of our house around the time that the alarm went off. Then common sense would suggest that the additional information concerning the occurrence of a small earthquake will somehow modify our probability assessment of there actually being a burglar in the house. 

We may use Bayesian probability theory to examine how the knowledge of a small earthquake having occurred translates to our state of knowledge regarding the plausibility of a burglary. The narrative we will formally analyze is taken from \cite{MacKay03}: 
\begin{quotation}
	\noindent Fred lives in Los Angeles and commutes 60 miles to work. Whilst at work, he receives a phone-call from his neighbor saying that Fred's burglar alarm is ringing. While driving home to investigate, Fred hears on the radio that there was a small earthquake that day near his home. 
\end{quotation}

The propositions that will go in our Bayesian inference network are:
\begin{align}
    B &= \text{Burglary,} \notag  \\
		\overline{B} &= \text{No burglary,} \notag  \\
		\notag  \\
		A &= \text{Alarm,} \notag  \\
		\overline{A} &= \text{No alarm,} \notag  \\
		\notag  \\
		E &= \text{Small earthquake,} \notag  \\
		\overline{E} &= \text{No earthquake,} \notag
\end{align}
where we will distinguish between two states of knowledge:
\begin{align}
    I_{1} &= \text{State of knowledge where hypothesis of earthquake is also entertained} \notag  \\
		I_{2}&= \text{State of knowledge where hypothesis of earthquake is not entertained} \notag  
\end{align}

We assume that the neighbor would never phone if the alarm is not ringing and that the radio report is trustworthy too; thus, we know for a fact that the alarm is ringing and that a small earthquake has occurred near the home. Furthermore, we assume that the occurrence of an earthquake and a burglary are independent. We also assume that a burglary alarm is almost certainly triggered by either a burglary or a small earthquake or both, that is,
\begin{equation}
	\label{eq.P1.C.6}
	P\!\left(\left.A\right|B \overline{E} I_{1}\right) = P\!\left(\left.A\right|\overline{B} E I_{1}\right) = P\!\left(\left.A\right|B E I_{1}\right) \rightarrow 1,
\end{equation}
whereas alarms in the absence of both a burglary and a small earthquake are extremely rare, that is,      
\begin{equation}
	\label{eq.P1.C.7}
	P\!\left(\left.A\right|\overline{B}\:\overline{E} I_{1}\right) \rightarrow 0.
\end{equation}
But if, in our state of knowledge, we do not entertain the possibility of an earthquake, then \eqref{eq.P1.C.6} and \eqref{eq.P1.C.7} will, respectively, collapse to
\begin{equation}
	\label{eq.P1.C.6bb}
	P\!\left(\left.A\right|B I_{2}\right) \rightarrow 1,
\end{equation}
and
\begin{equation}
	\label{eq.P1.C.7b}
	P\!\left(\left.A\right|\overline{B} I_{2}\right) \rightarrow 0.
\end{equation}

Let
\begin{equation}
	\label{eq.P1.C.8}
	P\!\left(E\right) = e, \qquad P\!\left(B\right) = b.
\end{equation}
Then we have, by way of the sum rule \eqref{eq.sum},
\begin{equation}
	\label{eq.P1.C.9}
	P\!\left(\overline{E}\right) = 1 - e, \qquad P\!\left(\overline{B}\right) = 1 - b.
\end{equation}

If we are in a state of knowledge where we allow for an earthquake, we have, by way of the product rule \eqref{eq.product}, as well as \eqref{eq.P1.C.6}, \eqref{eq.P1.C.7}, \eqref{eq.P1.C.8}, and \eqref{eq.P1.C.9}, that
\begin{align}
	\label{eq.P1.C.10}
	P\!\left(\left.A B \overline{E} \right| I_{1}\right) &= P\!\left(\left.A\right|B\overline{E} I_{1}\right) P\!\left(B\right) P\!\left(\overline{E}\right) \rightarrow b \left(1-e\right), \notag \\
	\notag \\
	P\!\left(\left.A \overline{B} E \right| I_{1}\right) &= P\!\left(\left.A\right|\overline{B} E I_{1}\right) P\!\left(\overline{B}\right) P\!\left(E\right) \rightarrow \left(1-b\right) e, \notag \\
	 \\
	P\!\left(\left.A B E \right| I_{1}\right) &= P\!\left(\left.A\right|B E I_{1}\right) P\!\left(B\right) P\!\left(E\right) \rightarrow b e, \notag \\
	\notag \\
	P\!\left(\left.A \overline{B}\:\overline{E}\right| I_{1} \right) &= P\!\left(\left.A\right|\overline{B}\:\overline{E} I_{1}\right) P\!\left(\overline{B}\right) P\!\left(\overline{E}\right) \rightarrow 0. \notag
\end{align}
By way of `marginalization', that is, an application of the generalized sum rule, \eqref{eq.genSum}, we obtain the probabilities
\begin{align}
	\label{eq.P1.C.11b}
	P\!\left(\left.A\overline{B} \right| I_{1}\right) &= P\!\left(\left.A \overline{B}\:\overline{E}\right| I_{1}\right) + P\!\left(\left.A \overline{B} E\right| I_{1}\right) \rightarrow \left(1-b\right) e, \notag \\
	\\
	P\!\left(\left.A B\right| I_{1}\right) &= P\!\left(\left.A B E\right| I_{1}\right) + 	P\!\left(\left.A B \overline{E}\right|I_{1}\right) \rightarrow b, \notag 
\end{align}
and
\begin{equation}	
 \label{eq.P1.C.11c}
P\!\left(\left.A\right|I_{1}\right) =  P\!\left(\left.A B\right|I_{1}\right) + P\!\left(\left.A\overline{B}\right|I_{1}\right) \rightarrow b + e - be,\end{equation}	
and
\begin{align}
	\label{eq.P1.C.11}
	P\!\left(\left.A\overline{E}\right|I_{1}\right) &= P\!\left(\left.A B \overline{E}\right|I_{1}\right) + P\!\left(\left.A \overline{B}\:\overline{E}\right|I_{1}\right) \rightarrow b \left(1-e\right), \notag \\
	\\
	P\!\left(\left.A E\right|I_{1}\right) &= P\!\left(\left.A B E\right|I_{1}\right) + P\!\left(\left.A \overline{B} E\right|I_{1}\right) \rightarrow e. \notag 
\end{align}

But if we are in a state of knowledge where we do not allow for an earthquake, we have, by way of the product rule \eqref{eq.product}, as well as \eqref{eq.P1.C.6bb},  \eqref{eq.P1.C.7b}, \eqref{eq.P1.C.8}, and \eqref{eq.P1.C.9}, that
\begin{align}
	\label{eq.P1.C.10b}
	P\!\left(\left.A B\right| I_{2}\right) &= P\!\left(\left.A\right|B I_{2}\right) P\!\left(B\right)  \rightarrow b , \notag \\
				\\
	P\!\left(\left.A \overline{B}\right|I_{2}\right) &= P\!\left(\left.A\right|\overline{B} I_{2}\right) P\!\left(\overline{B}\right) \rightarrow 0, \notag
\end{align}
By way of `marginalization', that is, an application of the generalized sum rule, \eqref{eq.genSum}, we obtain the probability
\begin{equation}
	\label{eq.P1.C.11bb}
	P\!\left(\left.A \right|I_{2}\right) = P\!\left(\left.A B \right| I_{2}\right) + 	P\!\left(\left.A \overline{B}\right| I_{2}\right) \rightarrow b. 
\end{equation}

The moment Fred hears that his burglary alarm is going off, then there are two possibilities. 

One possibility is that Fred may be new to Los Angeles and, consequently, overlook the possibility of a small earthquake triggering his burglary alarm, that is, his state of knowledge is $I_{2}$, which will make his prior probability of his alarm going off, go to \eqref{eq.P1.C.11bb}.

Fred then assesses, by way of the product rule \eqref{eq.product},\eqref{eq.P1.C.10b} and \eqref{eq.P1.C.11bb}, the likelihood of a burglary to be 
\begin{equation}
	\label{eq.P1.C.12b}
		P\!\left(\left.B\right|A I_{2}\right) = \frac{P\!\left(\left.A B\right|I_{2}\right)}{P\!\left(\left.A\right|I_{2}\right)} \rightarrow \frac{b}{b} = 1,
\end{equation}
which leaves him greatly distressed, as he drives to his home to investigate.

Another possibility is that Fred is a veteran Los Angeleno and, as a consequence, instantly will take into account the hypothesis of a small tremor occurring near his house, that is, his state of knowledge is $I_{1}$. 

Fred then assesses, by way of the product rule \eqref{eq.product},\eqref{eq.P1.C.11b} and \eqref{eq.P1.C.11c}, the likelihood of a burglary to be 
\begin{equation}
	\label{eq.P1.C.12c}
		P\!\left(\left.B\right|A I_{1}\right) = \frac{P\!\left(\left.A B\right|I_{1}\right)}{P\!\left(\left.A\right|I_{1}\right)} \rightarrow \frac{b}{b + e - b e} \approx  \frac{b}{b + e},
\end{equation}
seeing that $b + e >\,> b e$. 

If earthquakes are somewhat more common than burglaries, then Fred, based on his \eqref{eq.P1.C.12b}, may still hope for the best, as he drives home to investigate.

Either way, the moment that Fred hears on the radio that a small earthquake has occurred near his house, around the time when the burglary alarm went off, then, by way of the product rule \eqref{eq.product} and \eqref{eq.P1.C.10}  and \eqref{eq.P1.C.11}, Fred updates the likelihood of a burglary to be
\begin{equation}
	\label{eq.P1.C.12}
	P\!\left(\left.B\right|A E I_{1}\right) = \frac{P\!\left(\left.A B E\right|I_{1}\right)}{P\!\left(\left.A E\right|I_{1}\right)} \rightarrow \frac{b e}{e} = b.
\end{equation}
In the presence of an alternative explanation for the triggering of the burglary alarm, that is, a small earthquake occurring, the burglary alarm has lost its predictive power over the prior probability of a burglary, that is, \eqref{eq.P1.C.8} and \eqref{eq.P1.C.12},
\begin{equation}
	\label{eq.P1.C.13}
	P\!\left(\left.B\right|A E I_{1}\right) \rightarrow 	P\!\left(B\right).
\end{equation}

Consequently, Fred's fear for a burglary, as he rides home, after having heard that a small earthquake did occur, will only be dependent upon his assessment of the general likelihood of a burglary occurring. If we assume that Fred lives in a nice neighborhood, rather than some crime-ridden ghetto, then we can imagine that Fred will be, if not greatly, then at least somewhat, relieved. 

One of the arguments made against Bayesian probability theory as a normative model for human rationality is that people are generally numerical illiterate. Hence, the Bayesian model is deemed to be too numerical a model for human inference, \cite{Slovic04}. 

However, note that the Bayesian analysis given here was purely qualitative, in that no actual numerical values were given to our probabilities, apart from \eqref{eq.P1.C.6}, \eqref{eq.P1.C.6bb}, \eqref{eq.P1.C.7}, and  \eqref{eq.P1.C.7b}, which are limit cases of certainty and, hence, in a sense, may also be considered to be qualitative. 

Moreover, the result of this qualitative analysis seems to be intuitive enough. Indeed, the qualitative correspondence of the product and sum rules with common sense has been noted and demonstrated time and again by many distinguished scientists, including Laplace \cite{Laplace1819}, Keynes \cite{Keynes21}, Jeffreys \cite{Jeffreys39}, Polya \cite{Polya45, Polya54}, Cox \cite{Cox61}, Tribus \cite{Tribus69}, de~Finetti \cite{deFinetti74}, Rosenkrantz \cite{Rosenkrantz77}, and Jaynes \cite{Jaynes03}.

As an aside, Jaynes warns us in \cite{Jaynes03} that in any Bayesian inference problem we must conditionalize on our state of knowledge $I_{i}$. These authors were very well aware of this warning. Nonetheless, in our initial discussion of this toy problem we did heed Jaynes' warning. Because of the perceived gain in notational compactness. But to ignore Jaynes' warnings is to invite the potential for confusion, as we ourselves found out.

This object lesson notwithstanding, we will not conditionalize on our state of knowledge in our further discussion of the Bayesian decision theory. Our excuse for this lack in Bayesian rigor is that our decisions will be assumed to be based on an `universal' prior state of knowledge, \cite{Skilling05}, which will allow us, for the sake of notational compactness, to drop the extra symbol $I_{i}$.

\section{Bernoulli's Error}
\label{Error}
Kahneman dedicates in his Nobel lecture a whole section on `Bernoulli's error' and on how prospect theory may remedy this error, \cite{Kahneman02}. 

In Kahneman's Nobel lecture we may read the following on Bernoulli's error:
\begin{quote}
\noindent The idea that decision makers evaluate outcomes by the utility of final asset positions has been retained in economic analyses for almost 300 years. This is rather remarkable, because the idea is easily shown to be wrong; I call it Bernoulli's error. Bernoulli's model is flawed because it is reference-independent: it assumes that the value that is assigned to a given state of wealth does not vary with the decision maker's initial state of wealth. 
\end{quote}
So, Bernoulli's model is claimed to be in error in that it would evaluate outcomes by the utility of final asset positions alone, without taking into account the initial wealth of the decision maker. 

But it may be checked, Paragraph~10 of \cite{Bernoulli38}, that Bernoulli gives an utility function of the form:
\begin{equation}
	\label{eq.P1.8.7a1}
	u\left(\left.S\right|S_{0}\right) = q \log \frac{S}{S_{0}}.
\end{equation}
where, adopting the Kahneman's terminology, $S$ and $S_{0}$ are, respectively, the final and initial asset states. Let the asset increment $\Delta S$ be defined as 
\begin{equation}
	\label{eq.P1.8.7a2}
	\Delta S = S - S_{0}.
\end{equation}
Then, by substituting \eqref{eq.P1.8.7a2} into \eqref{eq.P1.8.7a1}, we obtain the equivalent utility function, \eqref{eq.P1.5.4b}:  
\begin{equation}
	\label{eq.P1.8.7a3}
	u\left(\left.\Delta S\right|S_{0}\right) = q \log \frac{S_{0} + \Delta S}{S_{0}}.
\end{equation}
It follows that in Bernoulli's expected utility theory asset increments $\Delta S$, be they positive or negative, are evaluated relative to the initial wealth $S_{0}$ of the decision maker. 

This then begs the question: What was it, that led Kahneman to the misguided belief that in Bernoulli's model the gains and losses are not the carriers of utility?

Kahneman and Tversky state the first two tenets\footnote{The third, and last, tenet is that of \textit{loss aversion}, which states that $u$ must be some concave down function. The Bernoulli utility function is concave down.} of expected utility theory to be, respectively, the tenets of \textit{expectation} and \textit{asset integration}, \cite{Kahneman79}. The tenet of expectation is
\begin{equation}
	\label{eq.P1.8.7a}
	U\!\left(x_{1}, p_{1}; \ldots; x_{n}, p_{n}\right) = p_{1} u\!\left(x_{1}\right) + \cdots + p_{n} u\!\left(x_{n}\right).
\end{equation}
The tenet of asset integration states that if $w$ is our current asset position, that is, our initial wealth, then we will accept an uncertain prospect having outcomes $x_{i}$ if 
\begin{equation}
	\label{eq.P1.8.7a4}
	U\!\left(w + x_{1}, p_{1}; \ldots; w + x_{n}, p_{n}\right) > u\!\left(w\right).
\end{equation}
By substituting \eqref{eq.P1.8.7a} into \eqref{eq.P1.8.7a4}, we obtain the implied asset integration tenet:
\begin{equation}
	\label{eq.P1.8.7b}
	p_{1} u\!\left(w + x_{1}\right) + \cdots + p_{n} u\!\left(w + x_{n}\right) > u\!\left(w\right).
\end{equation}
But then we have that the tenets of expectation and asset integration, as stated by Kahneman and Tversky, are incompatible with Bernoulli's expected utility theory. 

For the expectation tenet \eqref{eq.P1.8.7a}, Bernoulli's expected utility theory \cite{Bernoulli38} implies the function $u$:
\begin{equation}
	\label{eq.P1.8.7c}
	 u\!\left(x\right) = q \log \frac{w + x}{w},
\end{equation}
whereas, under the asset integration tenet \eqref{eq.P1.8.7b}, the implied function $u$ would be
\begin{equation}
	\label{eq.P1.8.7d}
	 u\!\left(x\right) = q \log x.
\end{equation}
By way of \eqref{eq.P1.8.7c}, \eqref{eq.P1.8.7d}, and the fact that  
\[
		q \log \frac{w + x}{w} \neq q \log x,
\]
it is then demonstrated that the two tenets, as proposed by Kahneman and Tversky in their \cite{Kahneman79}, are incompatible with Bernoulli's expected utility theory. 

By dropping the expectation tenet \eqref{eq.P1.8.7a}, while retaining the asset integration in its form of \eqref{eq.P1.8.7b}, we may, very easily, do away with the Kahneman and Tversky inconsistency. 

By substituting the implied \eqref{eq.P1.8.7d} into the asset integration tenet \eqref{eq.P1.8.7b}, we find\footnote{Note that Laplace \cite{Laplace1819} discusses Bernoulli's suggestion by way of the equivalent inequality:
\[
	p_{1} \log\!\left(w + x_{1}\right) + \cdots + p_{n} \log\!\left(w + x_{n}\right) > \log\!\left(w\right).
\]}:
\begin{equation}
	\label{eq.P1.8.7e}
	p_{1} \left[q \log\!\left(w + x_{1}\right)\right] + \cdots + p_{n} \left[q \log\!\left(w + x_{n}\right)\right] > q \log\!\left(w\right),
\end{equation}
Then, by way of the properties of the log function, we may rewrite \eqref{eq.P1.8.7e} into the equivalent
\begin{equation}
	\label{eq.P1.8.7f}
	p_{1} \left[q \log \frac{w + x_{1}}{w}\right] + \cdots + p_{n} \left[q \log \frac{w + x_{n}}{w}\right] > 0,
\end{equation}
which for any psychologist should be recognizable as the weighted sum of Weber-Fechner utilities\footnote{The Weber-Fechner law is used, amongst other things, to determine the decibel scale of human sound perception, where the Weber constant has been experimentally determined as $q = \frac{10}{\log 10}$.}. 

Especially, if those psychologists, like Kahneman and Tversky, explicitly state that the facts of perceptual adaptation were in their minds when they began their joint research on decision making under risk, \cite{Kahneman02}.

So, if it is claimed by Kahneman and Tversky \cite{Kahneman79} that Bernoulli's model is in error, as it would evaluate outcomes by the utility of final asset states alone, rather than gains or losses. Then we may infer that Kahneman and Tversky did not fully realize\footnote{In all fairness, we ourselves initially thought that we had improved on the insurance example given by Jaynes in his \cite{Jaynes03}, by using the Weber-Fechner law, which we still remembered from our psychology days, rather than Laplace's
\[
	p_{1} \left[\log\!\left(w + x_{1}\right)\right] + \cdots + p_{n} \left[\log\!\left(w + x_{n}\right)\right] > \log\!\left(w\right).
\]
But it was professor Han Vrijling, to whom we owe a debt of gratitude, who first pointed us to Bernoulli's \cite{Bernoulli38}, and the equivalence of the Weber-Fechner law and the Bernoulli utility function. It was only then, that we realized that Laplace's formulation is equivalent to \eqref{eq.P1.8.7f}.} that, according to Bernoulli \cite{Bernoulli38}, their abstract \eqref{eq.P1.8.7b} necessarily implies the concrete \eqref{eq.P1.8.7f}. 

At the end of the quote on Bernoulli's error\footnote{Given at the beginning of the chapter.} we may find the following cryptic footnote by Kahneman:
\begin{quote}
\noindent What varies with wealth in Bernoulli's theory is the response to a given \textit{change} of wealth. This variation is represented by the curvature of the utility function for wealth. Such a function cannot be drawn if the utility of wealth is reference-dependent, because utility then depends not only on current wealth but also on the reference level of wealth.
\end{quote}
We now will offer up an interpretation of what is stated in this footnote; as it may shed some further light on the Kahneman and Tversky position.

Let us assume that Kahneman and Tversky had at least some sense of what is written in, and here we quote Kahneman, \cite{Kahneman02}, `the brilliant essay that introduced the first version of expected utility theory (Bernoulli, 1738)'; that is, we assume that they were aware of the fact that Bernoulli proposes to use the log function, in some shape or form, in order to assign utilities to outcomes. 

Then it may well have been that Kahneman and Tversky were under the wrongful impression that Benoulli's utility function is given as 
\begin{equation}
\label{eq.P1.8.7g}
	u = q \log\left(w + x\right),
\end{equation}

The first two sentences in the footnote then may be interpreted as expressing the idea that, with differing levels of wealth $w$, the supposed utility function \eqref{eq.P1.8.7g} will be more or less linear in a given change of wealth $x$. 

If in the third sentence we let `current wealth' stand for change in wealth and `reference level of wealth' for initial state of wealth, then we may interpret it as stating that \eqref{eq.P1.8.7g} misses the necessary structure to take into account the initial state of wealth in its utility assignments.

Bernoulli's supposed error then would be that he proposed as his utility function \eqref{eq.P1.8.7g}. But Bernoulli proposed \eqref{eq.P1.8.7c} instead of \eqref{eq.P1.8.7g}. 

The erroneous utility function \eqref{eq.P1.8.7g} is problematic in that it cannot assign a value of zero to a change of wealth of $x = 0$. In Kahneman and Tversky's prospect theory we have that changes in wealth $x$ are assigned values by way of the value function $v$, where\footnote{Note, Kahneman and Tversky's two-part value function is given as, \cite{Tversky92}:
\[
	v\!\left(x\right) = \begin{cases} x^{\alpha}, \qquad &x \geq 0 \\
																	-\lambda\left(-x\right)^{\beta}, \qquad &x < 0
											\end{cases} 
\].}
\begin{equation}
\label{eq.P1.8.7h}
 		v\!\left(x\right) = 0,				\qquad \text{for } x = 0.
\end{equation}
So, if we read in \cite{Kahneman02}:
\begin{quote}
\noindent Preferences appeared to be determined by attitudes to gains and losses, defined relative to a reference point, but Bernoulli's theory and its successors did not incorporate a reference point. We therefore proposed an alternative theory of risk, in which the carriers of utility are gains and losses - changes of wealth rather than states of wealth. Prospect theory (Kahneman \& Tversky, 1979) embraces the idea that preferences are reference-dependent, and includes the extra parameter that is required by this assumption.
\end{quote}
Then we may interpret is as saying that the erroneous $u$ cannot assign zero utilities to zero outcomes, whereas Kahneman and Tversky's value function $v$ can. 

It would then follow that that which is embraced by prospect theory is the \textit{constraint}\footnote{Constraints are not the same as parameters.} \eqref{eq.P1.8.7h} on the value function $v$. But this constraint also holds, trivially, for Bernoulli's utility function \eqref{eq.P1.8.7c}. 

Now, after having established some tentative understanding into the reasoning process that might have led Kahneman and Tversky to their misunderstanding of Bernoulli's position, and after having provided a possible interpretation of Kahneman's footnote, we may start to wonder: We know how the initial wealth $w$ factors into Bernoulli's expected utility theory, but how does this initial wealth factor into prospect theory?

We quote Kahneman and Tversky \cite{Kahneman79}:
\begin{quote}   
\noindent The emphasis on changes as the carriers of value should not be taken to imply that the value of a particular change is independent of initial position. Strictly speaking, value should be treated as a function in two arguments: the asset position that serves as a reference point, and the magnitude of the change (positive or negative) from that reference point.
\end{quote}
And we could not agree more, though we ourselves would have dropped the `strictly speaking' modifier, as it weakens the imperative. Kahneman and Tversky continue: 
\begin{quote}   
\noindent However, the preference order of prospects is not greatly altered by small or even moderate variations in asset position. \ldots Consequently, the representation of value as a function in one argument generally provides a satisfactory approximation.
\end{quote}

So, it is Kahneman himself rather than Bernoulli, who does not take explicitly into account the decision maker's initial state of wealth\footnote{As may be checked in the previous footnote.}.

\section{The Case Against Bayes; A Reprensentative Example.}
\label{Represent}
The psychological paradigm of heuristics and biases originated as a reaction to the shortcomings of the mathematical expected utility theory. 

In the 1950's it was found that expected utility theory, the then dominant decision theory, failed to adequately model human decision making in certain instances, leading to such paradoxes as the Ellsberg and Allais paradox. Consequently, Edwards and his research team of PhD-students and post-docs took it upon themselves to remedy the situation, \cite{Phillips06}. 

Kahneman and Tversky, both post-docs under Edwards, proposed to construct a systematic theory about the psychology of uncertainty and judgment. In this psychological theory a handful of heuristics would replace the mathematical laws of chance as a more realistic model for subjective judgment of uncertainty. But what started as a parsimonious theory of human inference, consisting of only three heuristics and their associated biases, \cite{Kahneman02}, has proliferated into 20 heuristics and an impressive 170+ biases\footnote{Source Wikipedia, search `heuristics' and `list of cognitive biases'.}. 

Heuristics are said to be mental short cuts or `rules of thumb' humans use to do inference. It is theorized that, as we do not always have the time or resources to compare all the information at hand, we use heuristics to do inference quickly and efficiently. 

However, or so we are warned, even though these mental short cuts will be helpful most of the times, if used carelessly heuristics may lead to heuristic-induced biases, that is, systematic errors in reasoning, \cite{Kahneman02}. 

For example, when we use the representativeness heuristic then we estimate the likelihood of an event by comparing it to an existing prototype that already exists in our minds, \cite{Kahneman73}. 

Our prototype is what we think is the most relevant or typical example of a particular event or object. The bias associated with the representativeness heuristic is that when making judgments based on representativeness we are likely to overestimate the likelihood that the representative event will occur; just because an event or object is representative does not mean that it is more likely to occur.

In order to demonstrate this bias Kahneman and Tversky performed the following experiment, \cite{Kahneman73}. They divided the participants in their study up into three groups. 

The \textit{base-rate group} was asked to guess the percentage of all first-year graduate students in the following nine fields of specialization: business administration, computer science, engineering, humanities and education, law, library science, medicine, physics, and social sciences. 

The base-rate group estimated the highest percentage of graduate students, with 20\%, to be in humanities and education, and the second lowest percentage, with 7\%, to be in computer sciences.

The \textit{similarity group} was presented with the following profile:
\begin{quotation}
\noindent Tom W. is of high intelligence, although lacking in true creativity. He has a need for order and clarity, and for neat and tidy systems in which every detail finds its appropriate place. His writing is rather dull and mechanical, occasionally enlivened by somewhat corny puns and by flashes of imagination of the sci-fi type. He has a strong drive for competence. He seems to feel little sympathy for other people and does not enjoy interacting with others. Self-centered, he nonetheless has a deep moral sense.	
\end{quotation}
After which they were asked to rate how similar Tom was perceived to be to the typical graduate student in each of the nine graduate specializations. 

The similarity group assigned computer science the highest ranking position, with a mean similarity of 2.1, whereas humanities and education was assigned the second lowest ranking position, with a mean similarity ranking of 7.2. 

The \textit{prediction group} was given the same personality sketch of Tom as the similarity group, with the following additional information:
\begin{quotation}
\noindent The preceding personality sketch of Tom W. was written during Tom's senior year in high school by a psychologist, on the basis of projective tests. Tom W. is currently a graduate student.
\end{quotation}
Then they were asked to rank the nine fields of graduate specialization in order of the likelihood that Tom was now a graduate student in each of these fields. 

It was found that 95\% of the prediction group judged that Tom was more likely to study computer science than humanities and education. 

Since the likelihood rankings of the prediction group closely follow the similarity rankings of the similarity group, whereas they do not follow the base-rate estimates of the base-rate group, Kahneman and Tversky conclude that the representativeness heuristic must have been used by the participants in the prediction group, \cite{Kahneman73}.

However, the use of representativeness does not necessarily imply the representativeness heuristic. We quote Kahneman and Tversky on the representativeness heuristic, \cite{Kahneman72}:
\begin{quotation}
\noindent Our thesis is that, in many situations, an event $A_{1}$ is judged more probable than an event $A_{2}$ whenever $A_{1}$ appears more representative than $A_{2}$. In other words, the ordering of events by their subjective probabilities coincides with their ordering by representativeness.	
\end{quotation}
We now will give a formal analysis of the representativeness heuristic, as proposed by Kahneman and Tversky. 

Let
\begin{align}
    A_{1} &= \text{Computer Science Student} \notag  \\
	A_{2} &= \text{Humanities and Education Student} \notag  \\
	B &= \text{Profile} \notag
\end{align}
If Tom's psychological profile, $B$, is more representative of computer science students, $A_{1}$, than of humanities students, $A_{2}$, then
\begin{equation}
	\label{eq.P1.C.1}
	P\!\left(\left.B\right|A_{1}\right) > P\!\left(\left.B\right|A_{2}\right),
\end{equation}
or, equivalently,
\begin{equation}
	\label{eq.P1.C.2}
	\frac{P\!\left(\left.B\right|A_{1}\right)}{P\!\left(\left.B\right|A_{2}\right)} > 1.
\end{equation}
So, the representativeness heuristic is build upon the thesis:
\begin{equation}
	\label{eq.P1.C.3}
	\frac{P\!\left(\left.B\right|A_{1}\right)}{ P\!\left(\left.B\right|A_{2}\right)} > 1 \implies \frac{P\!\left(\left.A_{1}\right|B\right)}{P\!\left(\left.A_{2}\right|B\right)} > 1,
\end{equation}
where `$\implies$' is the symbol for logical implication.

However, thesis \eqref{eq.P1.C.3} is unfounded in that it does not follow directly from the rules of plausible reasoning. Moreover, Kahneman and Tversky seem to intuit as much; seeing that they warn us for the bias of base rate neglect, when using their representativeness thesis, \cite{Kahneman73}.

By taking the conclusion part in the thesis \eqref{eq.P1.C.3} as the starting point of a formal Bayesian analysis, we may find, by way of the product rule
\[
		P\!\left(A\right) P\!\left(\left.B\right|A\right) = P\!\left(A B\right) = P\!\left(B\right) P\!\left(\left.A\right|B\right),
\]
that
\begin{align}
	\label{eq.P1.C.4}
	\frac{P\!\left(\left.A_{1}\right|B\right)}{P\!\left(\left.A_{2}\right|B\right)} &= \frac{P\!\left(B\right)}{P\!\left(B\right)} \frac{P\!\left(\left.A_{1}\right|B\right)}{P\!\left(\left.A_{2}\right|B\right)}  \notag \\	
			\notag \\
			& =  \frac{P\!\left(A_{1} B\right)}{P\!\left(A_{2} B\right)} \\
		 \notag\\
		&=	\frac{P\!\left(A_{1}\right)}{P\!\left(A_{2}\right)} \frac{P\!\left(\left.B\right|A_{1}\right)}{P\!\left(\left.B\right|A_{2}\right)}. \notag 
\end{align}
It follows that the correct thesis, which actually does take into account the base rate, would be
\begin{equation}
	\label{eq.P1.C.5}
	 \frac{P\!\left(\left.B\right|A_{1}\right)}{ P\!\left(\left.B\right|A_{2}\right)} > \frac{P\!\left(A_{2}\right)}{P\!\left(A_{1}\right)} \implies \frac{P\!\left(\left.A_{1}\right|B\right)}{P\!\left(\left.A_{2}\right|B\right)} > 1.
\end{equation}

Kahneman and Tversky make in \cite{Kahneman73} the implicit assumption that the reported use of representativeness, that is, an evaluation and use of the odds \eqref{eq.P1.C.2}, as reported by the participants of the prediction group, necessarily implies their thesis \eqref{eq.P1.C.3}. This leads them to conclude that the participants must have used their representativeness heuristic. 

However, it can be seen that this assumption is incorrect, as the Bayesian thesis \eqref{eq.P1.C.5} also makes use of the odds in \eqref{eq.P1.C.2} and, thus, representativeness. Moreover, based on the reported use of representativeness and the experimental data, one could make the case that the participants in the experiment intuitively made use of the correct \eqref{eq.P1.C.5}, instead of the erroneous \eqref{eq.P1.C.3}. 

Kahneman and Tversky report that the prior odds for humanities and education against computer science were estimated by the participants to be, \cite{Kahneman73}:   
\begin{equation}
	\label{eq.P1.C.6a}
	\frac{P\!\left(A_{2}\right)}{P\!\left(A_{1}\right)} \approx 3. 
\end{equation}
Now, \eqref{eq.P1.C.4}, or, equivalently, thesis \eqref{eq.P1.C.5}, tells us that \eqref{eq.P1.C.6a} together with
\begin{equation}
	\label{eq.P1.C.6b}
	\frac{P\!\left(\left.B\right|A_{1}\right)}{P\!\left(\left.B\right|A_{2}\right)} > 3 
\end{equation}
implies the conclusion
\begin{equation}
	\label{eq.P1.C.7a}
	\frac{P\!\left(\left.A_{1}\right|B\right)}{P\!\left(\left.A_{2}\right|B\right)} > 1,
\end{equation}
or, equivalently,
\begin{equation}
	\label{eq.P1.C.8a}
	P\!\left(\left.A_{1}\right|B\right) > P\!\left(\left.A_{2}\right|B\right).
\end{equation}
So, if 95\% of the participants in the third group judged that Tom was more likely to study computer science than humanities, then we may infer that 95\% of the participants deemed the odds-inequality \eqref{eq.P1.C.6b} to hold, which in our opinion is not that far-fetched\footnote{Indeed, the ranking of computer sciences was, with a mean similarity of 2.1, more than three times higher than the ranking of humanities, which had a mean similarity ranking of only 7.2, \cite{Kahneman73}. Even though rankings do not translate easily to probabilities, as a \textit{sine qua non}, this ordering of similarity rankings still constitutes corroborating evidence for inequality \eqref{eq.P1.C.6b} to have held for the participants in the Kahneman and Tversky experiment.}. 

Kahneman and Tversky, however, by way of an informational accuracy argument, are of the opinion that the plausibility judgments of the participants in the third group `drastically violate the rules of the normative [that is, Bayesian] rules of prediction', seeing that the following considerations were ignored by the participants in the prediction group, \cite{Kahneman73}:
\begin{quote}
First given the notorious invalidity of projective personality tests, it is very likely that Tom W. was never in fact as compulsive and as aloof as his description suggests. Second, even if the description was valid when Tom W. was in high school, it may not longer be valid now that he is in graduate school. Finally, even if the description is still valid, there are probably more people who fit that description among students of humanities and education than among students of computer science, simply because there are so many more students in the former than in the latter field.  
\end{quote}

As to the last consideration, it follows from the Bayesian `heuristic' \eqref{eq.P1.C.5} that the inference of Tom being a computer science student implies the corollary inference that among all the graduate students who fit Tom's description there will be more students of computer science than students of humanities and education. Even if there are many more students in the field of humanities and education than in the field of computer science. 

The odds \eqref{eq.P1.C.6a} represent the ratio of humanities and education students to computer science students. Whereas the odds \eqref{eq.P1.C.6b} represent the ratio of computer students having a profile like Tom's to humanities and eduction students having a like profile. If the latter odds dominate the former, then we must conclude, by way of \eqref{eq.P1.C.4}, that there are probably more people who fit Tom's description among students of computer science than among students of humanities and education; that is, $P\!\left(A_{1} B\right) > P\!\left(A_{2} B\right)$. Which invalidates Kahneman's and Tversky's last consideration\footnote{And as a further consequence, we also grow to doubt Kahneman and Tversky's competency somewhat, when it comes to matters of the normative rules of prediction.}.

This, then, leaves us with the following arguments for the thesis that people tend to neglect the base rate, when taking the mental shortcut \eqref{eq.P1.C.3}, thus, violating the rules of the normative [that is, Bayesian] rules of prediction:
\begin{enumerate}
	\item given the notorious [that is, clinical] invalidity of projective personality tests, it is very likely that Tom W. was never in fact as compulsive and as aloof as his description suggests.
	\item even if the description was valid when Tom W. was in high school, it may not longer be valid now that he is in graduate school.  
\end{enumerate}

Now, it would seem that these arguments pertain to some other thesis, namely, that the participants of the prediction group should have disregarded the description of Tom altogether, as no data was actually given. But this alternative thesis, apart from it not being the issue\footnote{Though Kahneman and Tversky have made it the issue, by way of their informational accuracy argument.}, is not all that compelling.

Because, if, in answer to the first argument, we filter out those qualifications which might point to compulsiveness and aloofness, then we are left with the following profile for Tom:
\begin{itemize}
	\item high intelligence, 
	\item dull and mechanical writing, 
	\item lacking in true creativity, 
	\item corny puns, 
	\item flashes of imagination of the sci-fi type, 
	\item strong drive for competence, 
	\item deep moral sense,
\end{itemize}
which tells us quite a lot. 

It tells us that Tom is very bright, does not like to read, as he apparently is not that lyrical in his writing, is not very artistic, has a sense of humor, has a passion for sci-fi, is disciplined, and has a sense of justice\footnote{Moreover, what is left out also gives us some tentative information on Tom's psychologist.}. 

As to the second argument, which is also the hardest to answer. It is taught at the psychology courses, that personality traits tend to be stable over long periods of time. So, if Tom did not like to write in high school\footnote{A liking, admittedly, is not a personality trait, but still.}, than chances are that he would pick a field op specialization in which he would not have to read and write a lot. Which would make humanities and education less probable a field of choice, and computer sciences a more probable one.

But, lest we forget, the original thesis under discussion was that normative rules of prediction tend to be neglected, as people tend to neglect the base rate; not the alternative thesis that the prediction group should have disregarded the description of Tom, because of the clinical invalidity of projective personality tests and the possibility that Tom might have changed his personality over the course of the few years between high school and college. 

If Kahneman and Tversky wish to prove their initial thesis, then at this point, having presented their experimental data, they should proceed to demonstrate that their subjects could not possibly have used the normative, that is, Bayesian, rules of prediction, as those rules would have implied results other than those that were observed. Which they do not.

Having come to the end of our discussion of the representative heuristic, we find that the reported plausibility judgments by the prediction group are not inconsistent with a possible use of the Bayesian `heuristic; \eqref{eq.P1.C.6a} through \eqref{eq.P1.C.8a}, seeing that the odds \eqref{eq.P1.C.6b} may be assumed to lie in the realm of the possible. This leaves the Kahneman and Tversky experiment inconclusive.

\section{Non-Linear Preferences}
\label{Non-Linear Preferences}
Tversky and Kahneman \cite{Tversky92}, state that non-linear preferences constitute one of the minimal challenges that must be met by any adequate descriptive theory of choice. We shall explain.

If we have a bet which has the following outcome probability distributions
\begin{equation}
	\label{eq.P1.9.1}
	p\!\left(\left.O\right| D_{1}\right) = \begin{cases}  1, \qquad O = 1.000.000 
										   \end{cases}
\end{equation}
and
\begin{equation}
	\label{eq.P1.9.2}
	p\!\left(\left.O\right| D_{2}\right) = \begin{cases}  0.99, \qquad O = 5.000.000 \\
														  0.01, \qquad O = 0
										   \end{cases}
\end{equation}
then people will typically prefer the bet $D_{1}$ over $D_{2}$. 

In contrast, if we have a bet which has the following outcome probability distributions
\begin{equation}
	\label{eq.P1.9.3}
	 p\!\left(\left.O\right|D_{1}\right) = \begin{cases} 0.90, \qquad O = 1.000.000 \\
														 0.10, \qquad O = 0 
										   \end{cases}
\end{equation} 
and
\begin{equation}
	\label{eq.P1.9.4}
	 p\!\left(\left.O\right|D_{2}\right) = \begin{cases} 0.89, \qquad O = 5.000.000 \\
														 0.11, \qquad O = 0 
										   \end{cases}
\end{equation}
then people will typically prefer the bet $D_{2}$ over $D_{1}$. 

Allais gave this example, in a slightly altered form, to demonstrate that Savage's fifth axiom of independence does not hold, \cite{Allais53a}. 

According to Savage's independence axiom, which we most emphatically do not endorse, we may add for both \eqref{eq.P1.9.1} and \eqref{eq.P1.9.2} a probability mass of $0.10$ to the proposition $u = 0$, while subtracting that same probability mass for the respective propositions $u = 5.000.000$ and $u = 1.000.000$, leading to \eqref{eq.P1.9.3} and \eqref{eq.P1.9.4}, and still maintain the same problem of choice. But this assumption, as one would hope, is shown to be incorrect by the observed reversal in preferences from bet $D_{1}$ over $D_{2}$ to bet $D_{2}$ over $D_{1}$.

Now, according to Kahneman and Tversky the example by Allais not only refutes Savage's axiom of independence, but it also shows that the difference between probabilities of 0.99 and 1.00 has more impact on preferences than the difference between 0.10 and 0.11. 

Kahneman and Tversky find this observation to be so full of meaning that they deem it to be a psychological phenomenon in and of itself, and proceed to label it as the `certainty effect', \cite{Kahneman79}, which later turns into `non-linear preferences', \cite{Tversky92}. But for Bayesians the phenomenon of non-linear preferences is both not that special and not that new\footnote{The particular U-shape of the non-informative Jeffreys' prior for the parameter $\theta$ of the beta distribution, 
\[
	p\!\left(\theta\right) \propto \theta^{-1} \left(1 - \theta\right)^{-1},
\]
is a consequence of the non-linearity of probabilities. If we make a change of variable from $\theta$ to the log-odds $\omega = \log\left[ \theta / \left( 1 - \theta\right) \right]$, then it may be found that the implied non-informative prior of the log-odds $\omega$ is uniform, 
\[
	p\!\left(\omega\right) \propto \text{constant}.
\]
So, log-odds, having the whole infinity of the $x$-axis at their disposal, are linear; whereas probabilities, being forced to lie within the heavily constricted interval $\left[0, 1\right]$, are not.}. 

While working on the German enigma code, during World War II, Good\footnote{Good wrote about 2000 articles on Bayesian statistics, found throughout the statistical and philosophical literature starting in 1940. Workers in the field generally granted that every idea in modern statistics can be found expressed by him in one or more of these articles; but their sheer number made it impossible to find or cite them, and most are only one or two pages long, dashed off in an hour and never developed further. So, for many years, whatever one did in Bayesian statistics, one just conceded priority to Jack Good by default, without attempting the literature search, which would have required days. Finally, in 1983, a bibliography was provided of most of the first 1517 of these articles with a long index, so it is now possible to give proper acknowledgments of his works up to 1983, \cite{Jaynes03}.} and Turing introduced the `deciban' measure which is measured in decibans:
\begin{equation}
	\label{eq.P1.9.14}
	 \text{deciban}\!\left(P\right) = 10\log_{10} \frac{P}{1-P},
\end{equation}
and which gives the plausibility of a proposition being true, relative to it not being true, \cite{Good80}. 

For undecidedness, that is, for a fifty-fifty change of some hypothesis $A$ being true, we have
\begin{equation}
	\label{eq.P1.9.15}
	 P = 1 - P = 0.5.
\end{equation}
Substituting \eqref{eq.P1.9.15} in \eqref{eq.P1.9.14}, we find that undecidedness, \eqref{eq.P1.9.15}, corresponds with
\begin{equation}
	\label{eq.P1.9.16}
	 \text{deciban}\!\left(0.5\right) = 10\log_{10}\!\left(1\right) = 0.
\end{equation}
Just as 1 db sound represents the just noticeable difference relative to silence, so a $\pm 1$ deciban change in probability represents the just noticeable difference relative to undecidedness, \cite{Good50}. 

Using \eqref{eq.P1.9.14}, we find that the decibans associated with the probabilities 0.99 and 1.00 are, respectively,
\begin{equation}
	\label{eq.P1.9.17}
	 \text{deciban}\!\left(0.99\right) = 10\log_{10} \frac{0.99}{0.01} = 19.96
\end{equation}
and
\begin{equation}
	\label{eq.P1.9.18}
	 \text{deciban}\!\left(1.00\right) = 10\log_{10} \frac{1.00}{0.00} \rightarrow \infty.
\end{equation}
Now, \eqref{eq.P1.9.18} tells us that a probability 1.00 is a limit case of absolute certainty. And \eqref{eq.P1.9.17} tells us that a probability of 0.99 is not, representing just under 20 decibans of evidence for proposition $A$ being true. 

So, the difference in evidence for proposition $A$ being true for the probabilities 0.99 and 1.00 is much more than 1 deciban:
\begin{equation}
	\label{eq.P1.9.19}
	 \text{deciban}\!\left(1.00\right) - \text{deciban}\!\left(0.99\right) = \infty >> 1.  
\end{equation}
In contrast, the probabilities of 0.10 and 0.11 correspond with a less than 1 deciban difference in evidence, \eqref{eq.P1.9.14}: 	 					\begin{equation}
	\label{eq.P1.9.20}
	 \text{deciban}\!\left(0.11\right) - \text{deciban}\!\left(0.10\right) = 0.46 < 1.	 
\end{equation}

So, if we find that subjects prefer the second bet in the second collection of bets, \eqref{eq.P1.9.3} and \eqref{eq.P1.9.4}, then this also may be interpreted as follows. Subjects are indifferent to the difference in probabilities 0.10 and 0.11, as this difference represents a change of less then 1 deciban in the plausibility of hypothesis $A$ being true, \eqref{eq.P1.9.20}. So, all things being equal, subjects choose the bet with the highest potential payout of 5.000.000 dollars.

In closing, the deciban, \eqref{eq.P1.9.14}, represents the intuitive scale on which the plausibility of proposition $A$ being true, relative to it not being true, is judged; that is, the deciban is the scale of the numerically coded plausibilities, whereas the probability,
\begin{equation}
	\label{eq.P1.9.21}
	 P = \frac{10^{\frac{\text{deciban}\!\left(P\right)}{10}}}{1 + 10^{\frac{\text{deciban}\!\left(P\right)}{10}} },
\end{equation}
represents the non-intuitive `technical' scale, which follows from the quantification of our common sense, \cite{Knuth10, Skilling08}.

So, if the difference between probabilities of 0.99 and 1.00 has more impact on preferences than the difference between 0.10 and 0.11, as is found in psychological experiments, then this is reflective of the fact that the qualitative properties of the intuitive deciban-scale, up to a certain point, are retained in the transformation to the more technical probability-scale. 

We say up to a certain point, because probability theory, which makes use of the technical probabilities, is common sense amplified, having a much higher probability resolution than our human brains can ever hope to achieve. More concretely, we expect that for human probability perception the range of meaningful decibans is bounded somewhere around, say, $\pm 40$ deciban, \cite{Jaynes03}.

\section{The Framing Effect}
\label{Framing}
The assumption that preferences are not affected by variations of irrelevant features of options or outcomes is called invariance, \cite{Tversky86}. 

According to Kahneman and Tversky invariance is an essential aspect of rationality, which is violated in demonstrations of framing effects, \cite{Kahneman02}. Now, in order to discuss these framing effects, we will first have to discuss the topic of loss and gain adaptation\footnote{What we will call loss and gain adaptation, is called \textit{shifts of reference} by Kahneman and Tversky in their \cite{Kahneman79}.}. 

Imagine a business man who has lost 2000 in a business venture, and now is facing a choice between a sure gain of 1000 and an even chance to win 2000 or nothing, \cite{Kahneman79}. 

If this business man has not adapted to his loss, he is likely to add this loss to all the outcomes and, consequently, by way of the Weber-Fechner law, code the problem as a choice between the following utility distributions
\begin{equation}
	\label{eq.P1.10b.1}
	p\!\left(\left.u\right|D_{1}\right) = \begin{cases}
																						1.0, \qquad u~= q \log \frac{S_{0} - 1000}{S_{0}}
																				\end{cases}
\end{equation}
and
\begin{equation}
	\label{eq.P1.10b.2}
	p\!\left(\left.u\right|D_{2}\right) = \begin{cases}
																						0.5, \qquad u~= q \log \frac{S_{0} - 2000}{S_{0}} \\
																						0.5, \qquad u~= q \log \frac{S_{0} - 0}{S_{0}} 
																				\end{cases}
\end{equation}
where $S_{0}$ is the pre-loss asset position and $q$ is the Weber constant for money. 

Looking at the increments in assets, it is predicted that our business man will tend to prefer $D_{2}$ over $D_{1}$, as this is the preferred choice under risk seeking in the negative domain. It follows that a failure to adapt to losses may induce risk seeking in the negative domain, \cite{Kahneman79}. 

Stated differently, a person who has not made peace with his losses is likely to accept gambles that would be unacceptable to him otherwise. We may find support for this hypothesis by the observation that the tendency to bet on long shots will increase in the course of a betting day, \cite{McGlothlin56}. 

However, if our business man has adapted to his loss, then he will update his pre-loss asset position $S_{0}$ to an adjusted asset position $S_{0}^{\left(\text{adj.}\right)}$ in which the loss is discounted:
\[
	S_{0}^{\left(\text{adj.}\right)} = S_{0} - 2000
\] 
and code the problem as a choice between the utility distributions
\begin{equation}
	\label{eq.P1.10b.3}
	p\!\left(\left.u\right|D_{1}\right) = \begin{cases}
																						1.0, \qquad u~= q \log \frac{S_{0}^{\left(\text{adj.}\right)}+\ 1000}{S_{0}^{\left(\text{adj.}\right)}}
																				\end{cases}
\end{equation}
and
\begin{equation}
	\label{eq.P1.10b.4}
	p\!\left(\left.u\right|D_{2}\right) = \begin{cases}
																						0.5, \qquad u~= q \log \frac{S_{0}^{\left(\text{adj.}\right)}+ \ 0}{S_{0}^{\left(\text{adj.}\right)}}  \\
																						0.5, \qquad u~= q \log \frac{S_{0}^{\left(\text{adj.}\right)}+ \ 2000}{S_{0}^{\left(\text{adj.}\right)}}
																				\end{cases}
\end{equation}
where $S_{0}^{\left(\text{adj.}\right)}$ is the post-loss asset position and $q$ is the Weber constant for money. 

Again looking at the increments in assets, we see that the signs of these increments have reversed. By this reversal in the sign of the asset increments, we go from a risk seeking in the negative domain scenario to a risk aversion in the positive domain scenario, \cite{Kahneman79}. So, it is now predicted that our business man, having already adapted to his loss, will tend to reverse his preferences, and choose $D_{1}$ over $D_{2}$.

Having introduced the concepts of loss and gain adaptation and the adjusted initial wealth $S_{0}^{\left(\text{adj.}\right)}$, we now may turn to the discussion of the framing effect. 

Consider the following problems, which were presented to two different groups of subjects, \cite{Kahneman79}.
\begin{quotation}
\noindent \textbf{Group 1}:
In addition to whatever you own, you have been given 1000. You are now asked to choose between
\[
		p\!\left(\left.O\right|D_{1}\right) = \begin{cases}
																						0.5, \qquad O = 0  \\
																						0.5, \qquad O = 1000
																				\end{cases}
\]
and
\[
		p\!\left(\left.O\right|D_{2}\right) = \begin{cases}
																						1.0, \qquad O = 500
																				\end{cases}
\]
\end{quotation}

\begin{quotation}
\noindent \textbf{Group 2}:
In addition to whatever you own, you have been given 2000. You are now asked to choose between
\[
		p\!\left(\left.O\right|D_{1}\right) = \begin{cases}
																						0.5, \qquad O = -1000  \\
																						0.5, \qquad O = 0
																				\end{cases}
\]
and
\[
		p\!\left(\left.O\right|D_{2}\right) = \begin{cases}
																						1.0, \qquad O = -500
																				\end{cases}
\]
\end{quotation}

\noindent It is found that 84\% of $N = 70$ subjects in Group 1 prefer bet $D_{2}$ over bet $D_{1}$; whereas 69\% of $N = 68$ subjects in Group 2 prefer bet $D_{1}$ over bet $D_{2}$, indicating risk seeking in the negative domain, \cite{Kahneman79}. 

This result may indicate that the subjects in both groups adapted to the respective gifts of 1000 and 2000, by adjusting their initial $S_{0}$ into a new $S_{0}^{\left(\text{adj.}\right)}$, before choosing between the two bets. Since we would have expected to see the same preferences for both groups had the gifts been discounted in the outcomes\footnote{It may be checked that a discounting of the respective gifts in the corresponding outcomes would have resulted in identical outcome probability distributions for both Groups 1 and 2:
\[
		p\!\left(\left.O\right|D_{1}\right) = \begin{cases}
																						0.5, \qquad O = 1000  \\
																						0.5, \qquad O = 2000
																				\end{cases}
\]
and
\[
		p\!\left(\left.O\right|D_{2}\right) = \begin{cases}
																						1.0, \qquad O = 1500
																				\end{cases}.
\]}. But instead, we observe a reversal in preferences. So, we conclude that the gifts must have been discounted in initial wealth, rather than in the outcomes. 

Furthermore, based on the unadjusted outcomes alone, as given in the corresponding outcome probability distributions, we would expect both the observed risk aversion in the positive domain in Group 1, that is, the observed preference of $D_{2}$ over $D_{1}$, and the observed risk seeking in the negative domain in Group 2, that is, the observed preference of $D_{1}$ over $D_{2}$. This then also points to an adjustment of the initial wealth, rather an adjustment of the outcomes.

However, an alternative, more parsimonious, explanation for this framing effect, or, equivalently, the observed reversal in preferences in both groups, would be that the respective gifts of 1000 and 2000 were neglected by the subjects. And we would tend to agree with Kahneman and Tversky on this one, \cite{Kahneman79}. For, when reviewing these hypothetical choices, we ourselves overlooked these gifts too. 

But where Kahneman and Tversky see the framing effect as an indication that such monetary gifts, as a rule, will not factor into our real-life decisions, we quote, \cite{Kahneman79}:
\begin{quote}
The apparent neglect of a bonus [our gift] that was common to both options [our decisions $D_{1}$ and $D_{2}$] in Problems 11 and 12 [our Groups~1 and~2] implies that the carriers of value or utility are changes of wealth, rather than final asset positions that include current wealth. This conclusion is the cornerstone of an alternative theory of risky choice [their prospect theory].       
\end{quote}
We, instead, propose that this neglect of the gifts point to the limitations of the experimental method of hypothetical choices, as employed by Kahneman and Tversky.

Strictly speaking for ourselves, the receiving of a real-life gift of either 1000 or 2000 euros would be quite the occasion. Consequently, we can hardly imagine neglecting such a substantial sum of money, or, for that matter, not factoring its occurrence in our subsequent monetary decisions\footnote{The question if such a gift would be discounted into our initial wealth $S_{0}$, or into the specific outcomes $\Delta S$ of some set of decisions $D_{i}$, would be dependent on the particular context in which the gift was received.}. 

So, it would seem, at least based on the above \cite{Kahneman79} quotation, that Kahneman and Tversky build their prospect theory around a phenomenon which, as our introspection would suggest, is nothing but an experimental artifact of the method of hypothetical choices.

\end{document}